\documentclass[11pt,twocolumn]{article}
\usepackage{amsmath,amssymb,epic,eepic,rotating,color}
\usepackage[french]{babel}
\usepackage[T1]{fontenc}
\title{
\vspace{.5ex}
\begin{center}
\Huge \textbf{\textsf{Quelques mélanges parfaits de cartes}}
\end{center}\vspace{.21ex}}

\author{\textit{par} Aimé LACHAL
\footnote{
Adresse~: \hspace*{.35\textwidth} \mbox{}
\mbox{\textsc{Institut National des Sciences Appliquées de Lyon}}
\mbox{\textsl{Pôle de Mathématiques}}\hspace*{.15\textwidth} \mbox{}
\mbox{Bâtiment Léonard de Vinci, 20 avenue Albert Einstein}
\mbox{69621 Villeurbanne Cedex, \textsc{France}}
\mbox{E-mail: \tt{aime.lachal@insa-lyon.fr}}
}
}
\date{}

\newlength{\centrage}
\divide\centrage by 2
\oddsidemargin \centrage

\newtheorem{pr}{Proposition}
\newtheorem{theo}[pr]{Théorème}
\newtheorem{lm}[pr]{Lemme}
\newtheorem{co}[pr]{Corollaire}

\newcommand{\bpr}[1]{\begin{pr}#1\end{pr}}
\newcommand{\bth}[1]{\begin{theo}#1\end{theo}}

\newcommand{\bco}[1]{\begin{co}#1\end{co}}
\newcommand{\beq}{\begin{equation}}
\newcommand{\eeq}{\end{equation}}
\newcommand{\beqa}{\begin{eqnarray*}}
\newcommand{\eeqa}{\end{eqnarray*}}
\newcommand{\beqan}{\begin{eqnarray}}
\newcommand{\eeqan}{\end{eqnarray}}
\newcommand{\bitem}{\begin{itemize}}
\newcommand{\eitem}{\end{itemize}}

\newcommand{\dem}{\noindent {\sc Démonstration. }}
\newcommand{\fin}{\unitlength=.08em
                \hfill$\begin{picture}(7,6)
                \put(1,0){\line(0,1){6}}
                \put(7,0){\line(0,1){6}}
                \put(1,0){\line(1,0){6}}
                \put(1,6){\line(1,0){6}}
                \end{picture}$}
\newcommand{\rem}{\noindent {\sc Remarque. }}
\newcommand{\exam}{\noindent {\sc Exemple. }}

\newcommand{\e}{\varepsilon}
\newcommand{\f}{\varphi}
\newcommand{\ff}{\mathfrak{f}}
\newcommand{\tf}{\tilde{f}}
\newcommand{\tg}{\tilde{g}}
\renewcommand{\gg}{\mathfrak{g}}
\renewcommand{\th}{\tilde{h}}
\renewcommand{\ge}{\geqslant}
\renewcommand{\le}{\leqslant}
\newcommand{\N}{\ensuremath{\mathbb{N}}}
\newcommand{\Z}{\ensuremath{\mathbb{Z}}}
\newcommand{\cO}{\mathcal{O}}
\newcommand{\tcO}{\mathcal{\tilde{O}}}
\newcommand{\cS}{\mathcal{S}}
\newcommand{\card}{\mathop{\mathrm{card}}}
\newcommand{\modulo}[3]{\ensuremath{#1\equiv #2 \;[\hspace{-1em}\mod #3]}}
\newcommand{\modsimple}[2]{\ensuremath{#1\;{\hspace{-1em}\mod #2}}}
\newcommand{\modtext}[3]{\ensuremath{\mbox{$#1\equiv #2$}\;{[\hspace{-.65em}\mod #3]}}}

\newcommand{\dis}{\displaystyle}
\newcommand{\nos}{n$^{\mbox{\scriptsize os}}$}
\newcommand{\eg}{{e.g.}}
\newcommand{\ie}{{i.e.}}
\newcommand{\lqn}[1]{\noalign{\noindent $\displaystyle{#1}$}}
\newcommand{\scr}{\scriptscriptstyle}
\renewcommand{\text}{\textstyle}

\newcommand{\lb}{\linebreak}

\newcommand{\dessinjetonsun}
{\unitlength=.1\columnwidth
\begin{figure}[h!]
\begin{picture}(10,3.3)

\color{green}
\put(0.7,0.07){$\scr 1$}\put(1.57,0.07){$\scr V$}
\put(.5,0){\line(1,0){1.5}}\put(.5,.3){\line(1,0){1.5}}
\put(.5,0){\line(0,1){.3}}\put(2,0){\line(0,1){.3}}

\put(0.7,0.47){$\scr 2$}\put(1.57,0.47){$\scr V$}
\put(.5,.4){\line(1,0){1.5}}\put(.5,.7){\line(1,0){1.5}}
\put(.5,.4){\line(0,1){.3}}\put(2,.4){\line(0,1){.3}}

\put(0.7,0.87){$\scr 3$}\put(1.57,0.87){$\scr V$}
\put(.5,.8){\line(1,0){1.5}}\put(.5,1.1){\line(1,0){1.5}}
\put(.5,.8){\line(0,1){.3}}\put(2,.8){\line(0,1){.3}}

\put(0.7,1.27){$\scr 4$}\put(1.57,1.27){$\scr V$}
\put(.5,1.2){\line(1,0){1.5}}\put(.5,1.5){\line(1,0){1.5}}
\put(.5,1.2){\line(0,1){.3}}\put(2,1.2){\line(0,1){.3}}

\color{red}
\put(0.7,1.67){$\scr 5$}\put(1.57,1.67){$\scr R$}
\put(.5,1.6){\line(1,0){1.5}}\put(.5,1.9){\line(1,0){1.5}}
\put(.5,1.6){\line(0,1){.3}}\put(2,1.6){\line(0,1){.3}}

\put(0.7,2.07){$\scr 6$}\put(1.57,2.07){$\scr R$}
\put(.5,2){\line(1,0){1.5}}\put(.5,2.3){\line(1,0){1.5}}
\put(.5,2){\line(0,1){.3}}\put(2,2){\line(0,1){.3}}

\put(0.7,2.47){$\scr 7$}\put(1.57,2.47){$\scr R$}
\put(.5,2.4){\line(1,0){1.5}}\put(.5,2.7){\line(1,0){1.5}}
\put(.5,2.4){\line(0,1){.3}}\put(2,2.4){\line(0,1){.3}}

\put(0.7,2.87){$\scr 8$}\put(1.57,2.87){$\scr R$}
\put(.5,2.8){\line(1,0){1.5}}\put(.5,3.1){\line(1,0){1.5}}
\put(.5,2.8){\line(0,1){.3}}\put(2,2.8){\line(0,1){.3}}

\color{black}
\put(3.12,1.43){$\Longrightarrow$}

\color{green}
\put(5.2,0.97){$\scr 1$}\put(6.07,0.97){$\scr V$}
\put(5,.9){\line(1,0){1.5}}\put(5,1.2){\line(1,0){1.5}}
\put(5,.9){\line(0,1){.3}}\put(6.5,.9){\line(0,1){.3}}

\put(5.2,1.37){$\scr 2$}\put(6.07,1.37){$\scr V$}
\put(5,1.3){\line(1,0){1.5}}\put(5,1.6){\line(1,0){1.5}}
\put(5,1.3){\line(0,1){.3}}\put(6.5,1.3){\line(0,1){.3}}

\put(5.2,1.77){$\scr 3$}\put(6.07,1.77){$\scr V$}
\put(5,1.7){\line(1,0){1.5}}\put(5,2){\line(1,0){1.5}}
\put(5,1.7){\line(0,1){.3}}\put(6.5,1.7){\line(0,1){.3}}

\put(5.2,2.17){$\scr 4$}\put(6.07,2.17){$\scr V$}
\put(5,2.1){\line(1,0){1.5}}\put(5,2.4){\line(1,0){1.5}}
\put(5,2.1){\line(0,1){.3}}\put(6.5,2.1){\line(0,1){.3}}

\color{red}
\put(8.2,0.77){$\scr 5$}\put(9.07,0.77){$\scr R$}
\put(8,.7){\line(1,0){1.5}}\put(8,1){\line(1,0){1.5}}
\put(8,.7){\line(0,1){.3}}\put(9.5,.7){\line(0,1){.3}}

\put(8.2,1.17){$\scr 6$}\put(9.07,1.17){$\scr R$}
\put(8,1.1){\line(1,0){1.5}}\put(8,1.4){\line(1,0){1.5}}
\put(8,1.1){\line(0,1){.3}}\put(9.5,1.1){\line(0,1){.3}}

\put(8.2,1.57){$\scr 7$}\put(9.07,1.57){$\scr R$}
\put(8,1.5){\line(1,0){1.5}}\put(8,1.8){\line(1,0){1.5}}
\put(8,1.5){\line(0,1){.3}}\put(9.5,1.5){\line(0,1){.3}}

\put(8.2,1.97){$\scr 8$}\put(9.07,1.97){$\scr R$}
\put(8,1.9){\line(1,0){1.5}}\put(8,2.2){\line(1,0){1.5}}
\put(8,1.9){\line(0,1){.3}}\put(9.5,1.9){\line(0,1){.3}}

\color{black}
\put(7.06,1.43){$+$}

\end{picture}
\caption{Premier découpage}\label{fig-jetons1}
\end{figure}
}
\newcommand{\dessinjetonsdeux}
{\unitlength=.1\columnwidth
\begin{figure}[h!]
\begin{picture}(10,3.4)

\color{red}
\put(0.7,0.07){$\scr 5$}\put(1.57,0.07){$\scr R$}
\put(.5,0){\line(1,0){1.5}}\put(.5,.3){\line(1,0){1.5}}
\put(.5,0){\line(0,1){.3}}\put(2,0){\line(0,1){.3}}

\color{green}
\put(0.7,0.47){$\scr 1$}\put(1.57,0.47){$\scr V$}
\put(.5,.4){\line(1,0){1.5}}\put(.5,.7){\line(1,0){1.5}}
\put(.5,.4){\line(0,1){.3}}\put(2,.4){\line(0,1){.3}}

\color{red}
\put(0.7,0.87){$\scr 6$}\put(1.57,0.87){$\scr R$}
\put(.5,.8){\line(1,0){1.5}}\put(.5,1.1){\line(1,0){1.5}}
\put(.5,.8){\line(0,1){.3}}\put(2,.8){\line(0,1){.3}}

\color{green}
\put(0.7,1.27){$\scr 2$}\put(1.57,1.27){$\scr V$}
\put(.5,1.2){\line(1,0){1.5}}\put(.5,1.5){\line(1,0){1.5}}
\put(.5,1.2){\line(0,1){.3}}\put(2,1.2){\line(0,1){.3}}

\color{red}
\put(0.7,1.67){$\scr 7$}\put(1.57,1.67){$\scr R$}
\put(.5,1.6){\line(1,0){1.5}}\put(.5,1.9){\line(1,0){1.5}}
\put(.5,1.6){\line(0,1){.3}}\put(2,1.6){\line(0,1){.3}}

\color{green}
\put(0.7,2.07){$\scr 3$}\put(1.57,2.07){$\scr V$}
\put(.5,2){\line(1,0){1.5}}\put(.5,2.3){\line(1,0){1.5}}
\put(.5,2){\line(0,1){.3}}\put(2,2){\line(0,1){.3}}

\color{red}
\put(0.7,2.47){$\scr 8$}\put(1.57,2.47){$\scr R$}
\put(.5,2.4){\line(1,0){1.5}}\put(.5,2.7){\line(1,0){1.5}}
\put(.5,2.4){\line(0,1){.3}}\put(2,2.4){\line(0,1){.3}}

\color{green}
\put(0.7,2.87){$\scr 4$}\put(1.57,2.87){$\scr V$}
\put(.5,2.8){\line(1,0){1.5}}\put(.5,3.1){\line(1,0){1.5}}
\put(.5,2.8){\line(0,1){.3}}\put(2,2.8){\line(0,1){.3}}

\color{black}
\put(3.12,1.43){$\Longrightarrow$}

\color{red}
\put(5.2,0.97){$\scr 5$}\put(6.07,0.97){$\scr R$}
\put(5,.9){\line(1,0){1.5}}\put(5,1.2){\line(1,0){1.5}}
\put(5,.9){\line(0,1){.3}}\put(6.5,.9){\line(0,1){.3}}

\color{green}
\put(5.2,1.37){$\scr 1$}\put(6.07,1.37){$\scr V$}
\put(5,1.3){\line(1,0){1.5}}\put(5,1.6){\line(1,0){1.5}}
\put(5,1.3){\line(0,1){.3}}\put(6.5,1.3){\line(0,1){.3}}

\color{red}
\put(5.2,1.77){$\scr 6$}\put(6.07,1.77){$\scr R$}
\put(5,1.7){\line(1,0){1.5}}\put(5,2){\line(1,0){1.5}}
\put(5,1.7){\line(0,1){.3}}\put(6.5,1.7){\line(0,1){.3}}

\color{green}
\put(5.2,2.17){$\scr 2$}\put(6.07,2.17){$\scr V$}
\put(5,2.1){\line(1,0){1.5}}\put(5,2.4){\line(1,0){1.5}}
\put(5,2.1){\line(0,1){.3}}\put(6.5,2.1){\line(0,1){.3}}

\color{red}
\put(8.2,0.77){$\scr 7$}\put(9.07,0.77){$\scr R$}
\put(8,.7){\line(1,0){1.5}}\put(8,1){\line(1,0){1.5}}
\put(8,.7){\line(0,1){.3}}\put(9.5,.7){\line(0,1){.3}}

\color{green}
\put(8.2,1.17){$\scr 3$}\put(9.07,1.17){$\scr V$}
\put(8,1.1){\line(1,0){1.5}}\put(8,1.4){\line(1,0){1.5}}
\put(8,1.1){\line(0,1){.3}}\put(9.5,1.1){\line(0,1){.3}}

\color{red}
\put(8.2,1.57){$\scr 8$}\put(9.07,1.57){$\scr R$}
\put(8,1.5){\line(1,0){1.5}}\put(8,1.8){\line(1,0){1.5}}
\put(8,1.5){\line(0,1){.3}}\put(9.5,1.5){\line(0,1){.3}}

\color{green}
\put(8.2,1.97){$\scr 4$}\put(9.07,1.97){$\scr V$}
\put(8,1.9){\line(1,0){1.5}}\put(8,2.2){\line(1,0){1.5}}
\put(8,1.9){\line(0,1){.3}}\put(9.5,1.9){\line(0,1){.3}}

\color{black}
\put(7.06,1.43){$+$}

\end{picture}
\caption{Deuxième découpage}\label{fig-jetons2}
\end{figure}
}
\newcommand{\dessinjetonstrois}
{\unitlength=.1\columnwidth
\begin{figure}[h!]
\begin{picture}(10,3.4)

\color{red}
\put(0.7,0.07){$\scr 7$}\put(1.57,0.07){$\scr R$}
\put(.5,0){\line(1,0){1.5}}\put(.5,.3){\line(1,0){1.5}}
\put(.5,0){\line(0,1){.3}}\put(2,0){\line(0,1){.3}}

\put(0.7,0.47){$\scr 5$}\put(1.57,0.47){$\scr R$}
\put(.5,.4){\line(1,0){1.5}}\put(.5,.7){\line(1,0){1.5}}
\put(.5,.4){\line(0,1){.3}}\put(2,.4){\line(0,1){.3}}

\color{green}
\put(0.7,0.87){$\scr 3$}\put(1.57,0.87){$\scr V$}
\put(.5,.8){\line(1,0){1.5}}\put(.5,1.1){\line(1,0){1.5}}
\put(.5,.8){\line(0,1){.3}}\put(2,.8){\line(0,1){.3}}

\put(0.7,1.27){$\scr 1$}\put(1.57,1.27){$\scr V$}
\put(.5,1.2){\line(1,0){1.5}}\put(.5,1.5){\line(1,0){1.5}}
\put(.5,1.2){\line(0,1){.3}}\put(2,1.2){\line(0,1){.3}}

\color{red}
\put(0.7,1.67){$\scr 8$}\put(1.57,1.67){$\scr R$}
\put(.5,1.6){\line(1,0){1.5}}\put(.5,1.9){\line(1,0){1.5}}
\put(.5,1.6){\line(0,1){.3}}\put(2,1.6){\line(0,1){.3}}

\put(0.7,2.07){$\scr 6$}\put(1.57,2.07){$\scr R$}
\put(.5,2){\line(1,0){1.5}}\put(.5,2.3){\line(1,0){1.5}}
\put(.5,2){\line(0,1){.3}}\put(2,2){\line(0,1){.3}}

\color{green}
\put(0.7,2.47){$\scr 4$}\put(1.57,2.47){$\scr V$}
\put(.5,2.4){\line(1,0){1.5}}\put(.5,2.7){\line(1,0){1.5}}
\put(.5,2.4){\line(0,1){.3}}\put(2,2.4){\line(0,1){.3}}

\put(0.7,2.87){$\scr 2$}\put(1.57,2.87){$\scr V$}
\put(.5,2.8){\line(1,0){1.5}}\put(.5,3.1){\line(1,0){1.5}}
\put(.5,2.8){\line(0,1){.3}}\put(2,2.8){\line(0,1){.3}}

\color{black}
\put(3.12,1.43){$\Longrightarrow$}

\color{red}
\put(5.2,0.97){$\scr 7$}\put(6.07,0.97){$\scr R$}
\put(5,.9){\line(1,0){1.5}}\put(5,1.2){\line(1,0){1.5}}
\put(5,.9){\line(0,1){.3}}\put(6.5,.9){\line(0,1){.3}}

\put(5.2,1.37){$\scr 5$}\put(6.07,1.37){$\scr R$}
\put(5,1.3){\line(1,0){1.5}}\put(5,1.6){\line(1,0){1.5}}
\put(5,1.3){\line(0,1){.3}}\put(6.5,1.3){\line(0,1){.3}}

\color{green}
\put(5.2,1.77){$\scr 3$}\put(6.07,1.77){$\scr V$}
\put(5,1.7){\line(1,0){1.5}}\put(5,2){\line(1,0){1.5}}
\put(5,1.7){\line(0,1){.3}}\put(6.5,1.7){\line(0,1){.3}}

\put(5.2,2.17){$\scr 1$}\put(6.07,2.17){$\scr V$}
\put(5,2.1){\line(1,0){1.5}}\put(5,2.4){\line(1,0){1.5}}
\put(5,2.1){\line(0,1){.3}}\put(6.5,2.1){\line(0,1){.3}}

\color{red}
\put(8.2,0.77){$\scr 8$}\put(9.07,0.77){$\scr R$}
\put(8,.7){\line(1,0){1.5}}\put(8,1){\line(1,0){1.5}}
\put(8,.7){\line(0,1){.3}}\put(9.5,.7){\line(0,1){.3}}

\put(8.2,1.17){$\scr 6$}\put(9.07,1.17){$\scr R$}
\put(8,1.1){\line(1,0){1.5}}\put(8,1.4){\line(1,0){1.5}}
\put(8,1.1){\line(0,1){.3}}\put(9.5,1.1){\line(0,1){.3}}

\color{green}
\put(8.2,1.57){$\scr 4$}\put(9.07,1.57){$\scr V$}
\put(8,1.5){\line(1,0){1.5}}\put(8,1.8){\line(1,0){1.5}}
\put(8,1.5){\line(0,1){.3}}\put(9.5,1.5){\line(0,1){.3}}

\put(8.2,1.97){$\scr 2$}\put(9.07,1.97){$\scr V$}
\put(8,1.9){\line(1,0){1.5}}\put(8,2.2){\line(1,0){1.5}}
\put(8,1.9){\line(0,1){.3}}\put(9.5,1.9){\line(0,1){.3}}

\color{black}
\put(7.06,1.43){$+$}

\end{picture}
\caption{Troisième découpage}\label{fig-jetons3}
\end{figure}
}
\newcommand{\dessinjetonsquatre}
{\unitlength=.1\columnwidth
\begin{figure}[h!]
\begin{picture}(10,3.4)

\color{red}
\put(0.7,0.07){$\scr 7$}\put(1.57,0.07){$\scr R$}
\put(.5,0){\line(1,0){1.5}}\put(.5,.3){\line(1,0){1.5}}
\put(.5,0){\line(0,1){.3}}\put(2,0){\line(0,1){.3}}

\put(0.7,0.47){$\scr 8$}\put(1.57,0.47){$\scr R$}
\put(.5,.4){\line(1,0){1.5}}\put(.5,.7){\line(1,0){1.5}}
\put(.5,.4){\line(0,1){.3}}\put(2,.4){\line(0,1){.3}}

\put(0.7,0.87){$\scr 5$}\put(1.57,0.87){$\scr R$}
\put(.5,.8){\line(1,0){1.5}}\put(.5,1.1){\line(1,0){1.5}}
\put(.5,.8){\line(0,1){.3}}\put(2,.8){\line(0,1){.3}}

\put(0.7,1.27){$\scr 6$}\put(1.57,1.27){$\scr R$}
\put(.5,1.2){\line(1,0){1.5}}\put(.5,1.5){\line(1,0){1.5}}
\put(.5,1.2){\line(0,1){.3}}\put(2,1.2){\line(0,1){.3}}

\color{green}
\put(0.7,1.67){$\scr 3$}\put(1.57,1.67){$\scr V$}
\put(.5,1.6){\line(1,0){1.5}}\put(.5,1.9){\line(1,0){1.5}}
\put(.5,1.6){\line(0,1){.3}}\put(2,1.6){\line(0,1){.3}}

\put(0.7,2.07){$\scr 4$}\put(1.57,2.07){$\scr V$}
\put(.5,2){\line(1,0){1.5}}\put(.5,2.3){\line(1,0){1.5}}
\put(.5,2){\line(0,1){.3}}\put(2,2){\line(0,1){.3}}

\put(0.7,2.47){$\scr 1$}\put(1.57,2.47){$\scr V$}
\put(.5,2.4){\line(1,0){1.5}}\put(.5,2.7){\line(1,0){1.5}}
\put(.5,2.4){\line(0,1){.3}}\put(2,2.4){\line(0,1){.3}}

\put(0.7,2.87){$\scr 2$}\put(1.57,2.87){$\scr V$}
\put(.5,2.8){\line(1,0){1.5}}\put(.5,3.1){\line(1,0){1.5}}
\put(.5,2.8){\line(0,1){.3}}\put(2,2.8){\line(0,1){.3}}

\color{black}
\put(3.12,1.43){$\Longrightarrow$}

\color{red}
\put(5.2,0.97){$\scr 7$}\put(6.07,0.97){$\scr R$}
\put(5,.9){\line(1,0){1.5}}\put(5,1.2){\line(1,0){1.5}}
\put(5,.9){\line(0,1){.3}}\put(6.5,.9){\line(0,1){.3}}

\put(5.2,1.37){$\scr 8$}\put(6.07,1.37){$\scr R$}
\put(5,1.3){\line(1,0){1.5}}\put(5,1.6){\line(1,0){1.5}}
\put(5,1.3){\line(0,1){.3}}\put(6.5,1.3){\line(0,1){.3}}

\put(5.2,1.77){$\scr 5$}\put(6.07,1.77){$\scr R$}
\put(5,1.7){\line(1,0){1.5}}\put(5,2){\line(1,0){1.5}}
\put(5,1.7){\line(0,1){.3}}\put(6.5,1.7){\line(0,1){.3}}

\put(5.2,2.17){$\scr 6$}\put(6.07,2.17){$\scr R$}
\put(5,2.1){\line(1,0){1.5}}\put(5,2.4){\line(1,0){1.5}}
\put(5,2.1){\line(0,1){.3}}\put(6.5,2.1){\line(0,1){.3}}

\color{green}
\put(8.2,0.77){$\scr 3$}\put(9.07,0.77){$\scr V$}
\put(8,.7){\line(1,0){1.5}}\put(8,1){\line(1,0){1.5}}
\put(8,.7){\line(0,1){.3}}\put(9.5,.7){\line(0,1){.3}}

\put(8.2,1.17){$\scr 4$}\put(9.07,1.17){$\scr V$}
\put(8,1.1){\line(1,0){1.5}}\put(8,1.4){\line(1,0){1.5}}
\put(8,1.1){\line(0,1){.3}}\put(9.5,1.1){\line(0,1){.3}}

\put(8.2,1.57){$\scr 1$}\put(9.07,1.57){$\scr V$}
\put(8,1.5){\line(1,0){1.5}}\put(8,1.8){\line(1,0){1.5}}
\put(8,1.5){\line(0,1){.3}}\put(9.5,1.5){\line(0,1){.3}}

\put(8.2,1.97){$\scr 2$}\put(9.07,1.97){$\scr V$}
\put(8,1.9){\line(1,0){1.5}}\put(8,2.2){\line(1,0){1.5}}
\put(8,1.9){\line(0,1){.3}}\put(9.5,1.9){\line(0,1){.3}}

\color{black}
\put(7.06,1.43){$+$}

\end{picture}
\caption{Quatrième découpage}\label{fig-jetons4}
\end{figure}
}
\newcommand{\dessinjetonscinq}
{\unitlength=.1\columnwidth
\begin{figure}[h!]
\begin{picture}(10,3.4)

\color{green}
\put(0.7,0.07){$\scr 3$}\put(1.57,0.07){$\scr V$}
\put(.5,0){\line(1,0){1.5}}\put(.5,.3){\line(1,0){1.5}}
\put(.5,0){\line(0,1){.3}}\put(2,0){\line(0,1){.3}}

\color{red}
\put(0.7,0.47){$\scr 7$}\put(1.57,0.47){$\scr R$}
\put(.5,.4){\line(1,0){1.5}}\put(.5,.7){\line(1,0){1.5}}
\put(.5,.4){\line(0,1){.3}}\put(2,.4){\line(0,1){.3}}

\color{green}
\put(0.7,0.87){$\scr 4$}\put(1.57,0.87){$\scr V$}
\put(.5,.8){\line(1,0){1.5}}\put(.5,1.1){\line(1,0){1.5}}
\put(.5,.8){\line(0,1){.3}}\put(2,.8){\line(0,1){.3}}

\color{red}
\put(0.7,1.27){$\scr 8$}\put(1.57,1.27){$\scr R$}
\put(.5,1.2){\line(1,0){1.5}}\put(.5,1.5){\line(1,0){1.5}}
\put(.5,1.2){\line(0,1){.3}}\put(2,1.2){\line(0,1){.3}}

\color{green}
\put(0.7,1.67){$\scr 1$}\put(1.57,1.67){$\scr V$}
\put(.5,1.6){\line(1,0){1.5}}\put(.5,1.9){\line(1,0){1.5}}
\put(.5,1.6){\line(0,1){.3}}\put(2,1.6){\line(0,1){.3}}

\color{red}
\put(0.7,2.07){$\scr 5$}\put(1.57,2.07){$\scr R$}
\put(.5,2){\line(1,0){1.5}}\put(.5,2.3){\line(1,0){1.5}}
\put(.5,2){\line(0,1){.3}}\put(2,2){\line(0,1){.3}}

\color{green}
\put(0.7,2.47){$\scr 2$}\put(1.57,2.47){$\scr V$}
\put(.5,2.4){\line(1,0){1.5}}\put(.5,2.7){\line(1,0){1.5}}
\put(.5,2.4){\line(0,1){.3}}\put(2,2.4){\line(0,1){.3}}

\color{red}
\put(0.7,2.87){$\scr 6$}\put(1.57,2.87){$\scr R$}
\put(.5,2.8){\line(1,0){1.5}}\put(.5,3.1){\line(1,0){1.5}}
\put(.5,2.8){\line(0,1){.3}}\put(2,2.8){\line(0,1){.3}}

\color{black}
\put(3.12,1.43){$\Longrightarrow$}

\color{green}
\put(5.2,0.97){$\scr 3$}\put(6.07,0.97){$\scr V$}
\put(5,.9){\line(1,0){1.5}}\put(5,1.2){\line(1,0){1.5}}
\put(5,.9){\line(0,1){.3}}\put(6.5,.9){\line(0,1){.3}}

\color{red}
\put(5.2,1.37){$\scr 7$}\put(6.07,1.37){$\scr R$}
\put(5,1.3){\line(1,0){1.5}}\put(5,1.6){\line(1,0){1.5}}
\put(5,1.3){\line(0,1){.3}}\put(6.5,1.3){\line(0,1){.3}}

\color{green}
\put(5.2,1.77){$\scr 4$}\put(6.07,1.77){$\scr V$}
\put(5,1.7){\line(1,0){1.5}}\put(5,2){\line(1,0){1.5}}
\put(5,1.7){\line(0,1){.3}}\put(6.5,1.7){\line(0,1){.3}}

\color{red}
\put(5.2,2.17){$\scr 8$}\put(6.07,2.17){$\scr R$}
\put(5,2.1){\line(1,0){1.5}}\put(5,2.4){\line(1,0){1.5}}
\put(5,2.1){\line(0,1){.3}}\put(6.5,2.1){\line(0,1){.3}}

\color{green}
\put(8.2,0.77){$\scr 1$}\put(9.07,0.77){$\scr V$}
\put(8,.7){\line(1,0){1.5}}\put(8,1){\line(1,0){1.5}}
\put(8,.7){\line(0,1){.3}}\put(9.5,.7){\line(0,1){.3}}

\color{red}
\put(8.2,1.17){$\scr 5$}\put(9.07,1.17){$\scr R$}
\put(8,1.1){\line(1,0){1.5}}\put(8,1.4){\line(1,0){1.5}}
\put(8,1.1){\line(0,1){.3}}\put(9.5,1.1){\line(0,1){.3}}

\color{green}
\put(8.2,1.57){$\scr 2$}\put(9.07,1.57){$\scr V$}
\put(8,1.5){\line(1,0){1.5}}\put(8,1.8){\line(1,0){1.5}}
\put(8,1.5){\line(0,1){.3}}\put(9.5,1.5){\line(0,1){.3}}

\color{red}
\put(8.2,1.97){$\scr 6$}\put(9.07,1.97){$\scr R$}
\put(8,1.9){\line(1,0){1.5}}\put(8,2.2){\line(1,0){1.5}}
\put(8,1.9){\line(0,1){.3}}\put(9.5,1.9){\line(0,1){.3}}

\color{black}
\put(7.06,1.43){$+$}

\end{picture}
\caption{Cinquième découpage}\label{fig-jetons5}
\end{figure}
}
\newcommand{\dessinjetonssix}
{\unitlength=.1\columnwidth
\begin{figure}[h!]
\begin{picture}(10,3.4)

\color{green}
\put(0.7,0.07){$\scr 1$}\put(1.57,0.07){$\scr V$}
\put(.5,0){\line(1,0){1.5}}\put(.5,.3){\line(1,0){1.5}}
\put(.5,0){\line(0,1){.3}}\put(2,0){\line(0,1){.3}}

\put(0.7,0.47){$\scr 3$}\put(1.57,0.47){$\scr V$}
\put(.5,.4){\line(1,0){1.5}}\put(.5,.7){\line(1,0){1.5}}
\put(.5,.4){\line(0,1){.3}}\put(2,.4){\line(0,1){.3}}

\color{red}
\put(0.7,0.87){$\scr 5$}\put(1.57,0.87){$\scr R$}
\put(.5,.8){\line(1,0){1.5}}\put(.5,1.1){\line(1,0){1.5}}
\put(.5,.8){\line(0,1){.3}}\put(2,.8){\line(0,1){.3}}

\put(0.7,1.27){$\scr 7$}\put(1.57,1.27){$\scr R$}
\put(.5,1.2){\line(1,0){1.5}}\put(.5,1.5){\line(1,0){1.5}}
\put(.5,1.2){\line(0,1){.3}}\put(2,1.2){\line(0,1){.3}}

\color{green}
\put(0.7,1.67){$\scr 2$}\put(1.57,1.67){$\scr V$}
\put(.5,1.6){\line(1,0){1.5}}\put(.5,1.9){\line(1,0){1.5}}
\put(.5,1.6){\line(0,1){.3}}\put(2,1.6){\line(0,1){.3}}

\put(0.7,2.07){$\scr 4$}\put(1.57,2.07){$\scr V$}
\put(.5,2){\line(1,0){1.5}}\put(.5,2.3){\line(1,0){1.5}}
\put(.5,2){\line(0,1){.3}}\put(2,2){\line(0,1){.3}}

\color{red}
\put(0.7,2.47){$\scr 6$}\put(1.57,2.47){$\scr R$}
\put(.5,2.4){\line(1,0){1.5}}\put(.5,2.7){\line(1,0){1.5}}
\put(.5,2.4){\line(0,1){.3}}\put(2,2.4){\line(0,1){.3}}

\put(0.7,2.87){$\scr 8$}\put(1.57,2.87){$\scr R$}
\put(.5,2.8){\line(1,0){1.5}}\put(.5,3.1){\line(1,0){1.5}}
\put(.5,2.8){\line(0,1){.3}}\put(2,2.8){\line(0,1){.3}}

\color{black}
\put(3.12,1.43){$\Longrightarrow$}

\color{green}
\put(5.2,0.97){$\scr 1$}\put(6.07,0.97){$\scr V$}
\put(5,.9){\line(1,0){1.5}}\put(5,1.2){\line(1,0){1.5}}
\put(5,.9){\line(0,1){.3}}\put(6.5,.9){\line(0,1){.3}}

\put(5.2,1.37){$\scr 3$}\put(6.07,1.37){$\scr V$}
\put(5,1.3){\line(1,0){1.5}}\put(5,1.6){\line(1,0){1.5}}
\put(5,1.3){\line(0,1){.3}}\put(6.5,1.3){\line(0,1){.3}}

\color{red}
\put(5.2,1.77){$\scr 5$}\put(6.07,1.77){$\scr R$}
\put(5,1.7){\line(1,0){1.5}}\put(5,2){\line(1,0){1.5}}
\put(5,1.7){\line(0,1){.3}}\put(6.5,1.7){\line(0,1){.3}}

\put(5.2,2.17){$\scr 7$}\put(6.07,2.17){$\scr R$}
\put(5,2.1){\line(1,0){1.5}}\put(5,2.4){\line(1,0){1.5}}
\put(5,2.1){\line(0,1){.3}}\put(6.5,2.1){\line(0,1){.3}}

\color{green}
\put(8.2,0.77){$\scr 2$}\put(9.07,0.77){$\scr V$}
\put(8,.7){\line(1,0){1.5}}\put(8,1){\line(1,0){1.5}}
\put(8,.7){\line(0,1){.3}}\put(9.5,.7){\line(0,1){.3}}

\put(8.2,1.17){$\scr 4$}\put(9.07,1.17){$\scr V$}
\put(8,1.1){\line(1,0){1.5}}\put(8,1.4){\line(1,0){1.5}}
\put(8,1.1){\line(0,1){.3}}\put(9.5,1.1){\line(0,1){.3}}

\color{red}
\put(8.2,1.57){$\scr 6$}\put(9.07,1.57){$\scr R$}
\put(8,1.5){\line(1,0){1.5}}\put(8,1.8){\line(1,0){1.5}}
\put(8,1.5){\line(0,1){.3}}\put(9.5,1.5){\line(0,1){.3}}

\put(8.2,1.97){$\scr 8$}\put(9.07,1.97){$\scr R$}
\put(8,1.9){\line(1,0){1.5}}\put(8,2.2){\line(1,0){1.5}}
\put(8,1.9){\line(0,1){.3}}\put(9.5,1.9){\line(0,1){.3}}

\color{black}
\put(7.06,1.43){$+$}

\end{picture}
\caption{Sixième découpage}\label{fig-jetons6}
\end{figure}
}

\newcommand{\dessininshuffle}
{\unitlength=.1\columnwidth
\begin{figure}[h!]
\begin{picture}(10,3.8)
\put(0.7,0.18){$\scr 1$}\put(0.7,0.68){$\scr 2$}
\put(0.22,2.18){$\scr n-1$}\put(0.7,2.68){$\scr n$}
\put(3.61,-.08){$\scr n+1$}\put(3.61,0.42){$\scr n+2$}
\put(3.61,1.92){$\scr 2n-1$}\put(3.61,2.42){$\scr 2n$}

\put(0.54,3.5){\scriptsize 1\ier\ paquet}
\put(2.53,3.5){\scriptsize 2\ieme\ paquet}
\put(1,0.25){\line(1,0){1.5}}\put(1,0.75){\line(1,0){1.5}}
\put(1,2.25){\line(1,0){1.5}}\put(1,2.75){\line(1,0){1.5}}
\put(2,0){\line(1,0){1.5}}\put(2,0.5){\line(1,0){1.5}}
\put(2,2){\line(1,0){1.5}}\put(2,2.5){\line(1,0){1.5}}
\put(1.7,1.32){$\vdots$}\put(2.7,1.07){$\vdots$}

\put(4.8,1.3){$\Longrightarrow$}

\put(6.4,-.08){$\scr 1$}\put(6.4,.18){$\scr 2$}
\put(6.4,.42){$\scr 3$}\put(6.4,.68){$\scr 4$}
\put(5.75,1.92){$\scr 2n-3$}\put(5.75,2.18){$\scr 2n-2$}
\put(5.75,2.42){$\scr 2n-1$}\put(6.22,2.68){$\scr 2n$}

\put(8.33,-.08){$\scr n+1$}\put(8.33,.18){$\scr 1$}
\put(8.33,.42){$\scr n+2$}\put(8.33,.68){$\scr 2$}
\put(8.33,1.92){$\scr 2n-1$}\put(8.33,2.18){$\scr n-1$}
\put(8.33,2.43){$\scr 2n$}\put(8.33,2.68){$\scr n$}

\put(4.9,3.5){\scriptsize\no de place}\put(8.3,3.5){\scriptsize\no de carte}
\put(6.95,3.55){\vector(1,0){1}}\put(7.32,3.66){$\scr f$}
\put(6.7,0){\line(1,0){1.5}}\put(6.7,0.25){\line(1,0){1.5}}
\put(6.7,0.5){\line(1,0){1.5}}\put(6.7,0.75){\line(1,0){1.5}}
\put(6.7,2){\line(1,0){1.5}}\put(6.7,2.25){\line(1,0){1.5}}
\put(6.7,2.5){\line(1,0){1.5}}\put(6.7,2.75){\line(1,0){1.5}}
\put(7.4,1.19){$\vdots$}

\end{picture}
\caption{In-shuffle, permutation $f$}\label{fig-in-shuffle}
\end{figure}
}

\newcommand{\dessinoutshuffle}
{\unitlength=.1\columnwidth
\begin{figure}[h!]
\begin{picture}(10,3.8)

\put(0.7,-0.08){$\scr 1$}\put(0.7,0.43){$\scr 2$}
\put(0.22,1.97){$\scr n-1$}\put(0.7,2.43){$\scr n$}
\put(3.61,.18){$\scr n+1$}\put(3.61,0.68){$\scr n+2$}
\put(3.61,2.18){$\scr 2n-1$}\put(3.61,2.68){$\scr 2n$}

\put(0.54,3.5){\scriptsize 1\ier\ paquet}
\put(2.53,3.5){\scriptsize 2\ieme\ paquet}

\put(1,0){\line(1,0){1.5}}\put(1,0.5){\line(1,0){1.5}}
\put(1,2){\line(1,0){1.5}}\put(1,2.5){\line(1,0){1.5}}
\put(2,.25){\line(1,0){1.5}}\put(2,0.75){\line(1,0){1.5}}
\put(2,2.25){\line(1,0){1.5}}\put(2,2.75){\line(1,0){1.5}}
\put(1.7,1.07){$\vdots$}\put(2.7,1.32){$\vdots$}

\put(4.8,1.3){$\Longrightarrow$}

\put(6.4,-.08){$\scr 1$}\put(6.4,.18){$\scr 2$}
\put(6.4,.42){$\scr 3$}\put(6.4,.68){$\scr 4$}
\put(5.75,1.92){$\scr 2n-3$}\put(5.75,2.18){$\scr 2n-2$}
\put(5.75,2.42){$\scr 2n-1$}\put(6.22,2.68){$\scr 2n$}

\put(8.33,-.08){$\scr 1$}\put(8.33,.18){$\scr n+1$}
\put(8.33,.42){$\scr 2$}\put(8.33,.68){$\scr n+2$}
\put(8.33,1.92){$\scr n-1$}\put(8.33,2.18){$\scr 2n-1$}
\put(8.33,2.43){$\scr n$}\put(8.33,2.68){$\scr 2n$}

\put(4.9,3.5){\scriptsize\no de place}\put(8.3,3.5){\scriptsize\no de carte}
\put(6.95,3.55){\vector(1,0){1}}\put(7.32,3.66){$\scr g$}
\put(6.7,0){\line(1,0){1.5}}\put(6.7,0.25){\line(1,0){1.5}}
\put(6.7,0.5){\line(1,0){1.5}}\put(6.7,0.75){\line(1,0){1.5}}
\put(6.7,2){\line(1,0){1.5}}\put(6.7,2.25){\line(1,0){1.5}}
\put(6.7,2.5){\line(1,0){1.5}}\put(6.7,2.75){\line(1,0){1.5}}
\put(7.4,1.19){$\vdots$}

\end{picture}
\caption{Out-shuffle, permutation $g$}\label{fig-out-shuffle}
\end{figure}
}

\newcommand{\dessinmongeshuffleun}
{\unitlength=.1\columnwidth
\begin{figure}[h!]
\begin{picture}(10,3.8)

\put(0.75,3.5){\scriptsize paquet initial}
\put(0.7,-.08){$\scr 1$}\put(0.7,0.18){$\scr 2$}
\put(0.7,0.68){$\scr 4$}\put(0.7,.43){$\scr 3$}
\put(0.05,1.92){$\scr 2n-3$}\put(0.05,2.18){$\scr 2n-2$}
\put(0.05,2.42){$\scr 2n-1$}\put(0.52,2.68){$\scr 2n$}

\put(1,0){\line(1,0){1.5}}\put(1,.25){\line(1,0){1.5}}
\put(1,.5){\line(1,0){1.5}}\put(1,.75){\line(1,0){1.5}}
\put(1,2){\line(1,0){1.5}}\put(1,2.25){\line(1,0){1.5}}
\put(1,2.5){\line(1,0){1.5}}\put(1,2.75){\line(1,0){1.5}}
\put(1.7,1.18){$\vdots$}

\put(5.75,2.68){$\scr 2n-1$}\put(5.75,2.42){$\scr 2n-3$}
\put(6.4,.92){$\scr 4$}\put(6.4,1.43){$\scr 1$}
\put(6.4,1.18){$\scr 2$}\put(6.4,1.68){$\scr 3$}
\put(5.75,.18){$\scr 2n-2$}\put(6.22,-.08){$\scr 2n$}

\put(8.33,-.08){$\scr 1$}\put(8.33,.18){$\scr 2$}
\put(8.33,.92){$\scr n-1$}\put(8.33,1.18){$\scr n$}
\put(8.33,1.43){$\scr n+1$}\put(8.33,1.68){$\scr n+2$}
\put(8.33,2.43){$\scr 2n-1$}\put(8.33,2.68){$\scr 2n$}

\put(4.9,3.5){\scriptsize\no de carte}\put(8.3,3.5){\scriptsize\no de place}
\put(6.95,3.55){\vector(1,0){1}}\put(7.3,3.66){$\scr h_1$}
\put(6.7,0){\line(1,0){1.5}}\put(6.7,0.25){\line(1,0){1.5}}
\put(6.7,1){\line(1,0){1.5}}\put(6.7,1.25){\line(1,0){1.5}}
\put(6.7,1.5){\line(1,0){1.5}}\put(6.7,1.75){\line(1,0){1.5}}
\put(6.7,2.5){\line(1,0){1.5}}\put(6.7,2.75){\line(1,0){1.5}}
\put(7.4,.44){$\vdots$}\put(7.4,1.94){$\vdots$}

\color{red}
\put(3.3,.25){\line(2,1){2}}
\put(3.3,.75){\line(8,1){2}}
\put(3.3,2.25){\line(1,-1){2}}
\put(3.3,2.75){\line(8,-11){2}}
\dashline{.05}(2.65,.25)(3.15,.25)
\dashline{.05}(2.65,.75)(3.15,.75)
\dashline{.05}(2.65,2.25)(3.15,2.25)
\dashline{.05}(2.65,2.75)(3.15,2.75)

\dashline{.05}(5.45,0)(6.05,0)
\dashline{.05}(5.45,.25)(5.65,.25)
\dashline{.05}(5.45,1)(6.3,1)
\dashline{.05}(5.45,1.25)(6.3,1.25)

\color{blue}
\put(3.3,0){\line(4,3){2}}
\put(3.3,.5){\line(8,5){2}}
\put(3.3,2){\line(4,1){2}}
\put(3.3,2.5){\line(8,1){2}}
\dashline{.05}(2.65,0)(3.15,0)
\dashline{.05}(2.65,.5)(3.15,.5)
\dashline{.05}(2.65,2)(3.15,2)
\dashline{.05}(2.65,2.5)(3.15,2.5)

\dashline{.05}(5.45,1.5)(6.3,1.5)
\dashline{.05}(5.45,1.75)(6.3,1.75)
\dashline{.05}(5.45,2.5)(5.65,2.5)
\dashline{.05}(5.45,2.75)(5.65,2.75)

\end{picture}
\caption{Mélange de Monge, permutation $h_1$}\label{fig-Monge1}
\end{figure}
}
\newcommand{\dessinmongeshuffledeux}
{\unitlength=.1\columnwidth
\begin{figure}[h!]
\begin{picture}(10,3.8)

\put(0.75,3.5){\scriptsize paquet initial}
\put(0.7,-.08){$\scr 1$}\put(0.7,0.18){$\scr 2$}
\put(0.7,0.43){$\scr 3$}\put(0.7,.68){$\scr 4$}
\put(0.05,1.92){$\scr 2n-3$}\put(0.05,2.18){$\scr 2n-2$}
\put(0.05,2.42){$\scr 2n-1$}\put(0.52,2.68){$\scr 2n$}

\put(1,0){\line(1,0){1.5}}\put(1,.25){\line(1,0){1.5}}
\put(1,.5){\line(1,0){1.5}}\put(1,.75){\line(1,0){1.5}}
\put(1,2){\line(1,0){1.5}}\put(1,2.25){\line(1,0){1.5}}
\put(1,2.5){\line(1,0){1.5}}\put(1,2.75){\line(1,0){1.5}}
\put(1.7,1.18){$\vdots$}

\put(5.75,-.08){$\scr 2n-1$}\put(5.75,.18){$\scr 2n-3$}
\put(6.4,.92){$\scr 3$}\put(6.4,1.18){$\scr 1$}
\put(6.4,1.43){$\scr 2$}\put(6.4,1.68){$\scr 4$}
\put(5.75,2.42){$\scr 2n-2$}\put(6.22,2.68){$\scr 2n$}

\put(8.33,-.08){$\scr 1$}\put(8.33,.18){$\scr 2$}
\put(8.33,.92){$\scr n-1$}\put(8.33,1.18){$\scr n$}
\put(8.33,1.43){$\scr n+1$}\put(8.33,1.68){$\scr n+2$}
\put(8.33,2.43){$\scr 2n-1$}\put(8.33,2.68){$\scr 2n$}

\put(4.9,3.5){\scriptsize\no de carte}\put(8.3,3.5){\scriptsize\no de place}
\put(6.95,3.55){\vector(1,0){1}}\put(7.3,3.66){$\scr h_2$}
\put(6.7,0){\line(1,0){1.5}}\put(6.7,0.25){\line(1,0){1.5}}
\put(6.7,1){\line(1,0){1.5}}\put(6.7,1.25){\line(1,0){1.5}}
\put(6.7,1.5){\line(1,0){1.5}}\put(6.7,1.75){\line(1,0){1.5}}
\put(6.7,2.5){\line(1,0){1.5}}\put(6.7,2.75){\line(1,0){1.5}}
\put(7.4,.44){$\vdots$}\put(7.4,1.94){$\vdots$}

\color{red}
\put(3.3,0){\line(8,5){2}}
\put(3.3,.5){\line(4,1){2}}
\put(3.3,2){\line(8,-7){2}}
\put(3.3,2.5){\line(4,-5){2}}
\dashline{.05}(2.65,0)(3.15,0)
\dashline{.05}(2.65,.5)(3.15,.5)
\dashline{.05}(2.65,2)(3.15,2)
\dashline{.05}(2.65,2.5)(3.15,2.5)

\dashline{.05}(5.45,0)(5.65,0)
\dashline{.05}(5.45,.25)(5.65,.25)
\dashline{.05}(5.45,1)(6.3,1)
\dashline{.05}(5.45,1.25)(6.3,1.25)

\color{blue}
\put(3.3,.25){\line(8,5){2}}
\put(3.3,.75){\line(2,1){2}}
\put(3.3,2.25){\line(8,1){2}}
\put(3.3,2.75){\line(1,0){2}}
\dashline{.05}(2.65,.25)(3.15,.25)
\dashline{.05}(2.65,.75)(3.15,.75)
\dashline{.05}(2.65,2.25)(3.15,2.25)
\dashline{.05}(2.65,2.75)(3.15,2.75)

\dashline{.05}(5.45,1.5)(6.3,1.5)
\dashline{.05}(5.45,1.75)(6.3,1.75)
\dashline{.05}(5.45,2.5)(5.65,2.5)
\dashline{.05}(5.45,2.75)(6.05,2.75)

\end{picture}
\caption{Mélange de Monge, permutation $h_2$}\label{fig-Monge2}
\end{figure}
}

\newcommand{\dessinmongeshuffletrois}
{\unitlength=.1\columnwidth
\begin{figure}[h!]
\begin{picture}(10,3.8)

\put(0.7,0.18){$\scr 1$}\put(0.7,0.68){$\scr 2$}
\put(0.22,2.18){$\scr n-1$}\put(0.7,2.68){$\scr n$}
\put(3.61,-.08){$\scr 2n$}\put(3.61,0.42){$\scr 2n-1$}
\put(3.61,1.92){$\scr n+2$}\put(3.61,2.42){$\scr n+1$}

\put(0.54,3.5){\scriptsize 1\ier\ paquet}
\put(2.5,3.5){\scriptsize 2\ieme\ paquet}\put(2.6,3.2){\scriptsize retourné}
\put(1,0.25){\line(1,0){1.5}}\put(1,0.75){\line(1,0){1.5}}
\put(1,2.25){\line(1,0){1.5}}\put(1,2.75){\line(1,0){1.5}}
\put(2,0){\line(1,0){1.5}}\put(2,0.5){\line(1,0){1.5}}
\put(2,2){\line(1,0){1.5}}\put(2,2.5){\line(1,0){1.5}}
\put(1.7,1.32){$\vdots$}\put(2.7,1.07){$\vdots$}

\put(4.8,1.3){$\Longrightarrow$}

\put(6.4,-.08){$\scr 1$}\put(6.4,.18){$\scr 2$}
\put(6.4,.42){$\scr 3$}\put(6.4,.68){$\scr 4$}
\put(5.75,1.92){$\scr 2n-3$}\put(5.75,2.18){$\scr 2n-2$}
\put(5.75,2.42){$\scr 2n-1$}\put(6.22,2.68){$\scr 2n$}

\put(8.33,-.08){$\scr 2n$}\put(8.33,.18){$\scr 1$}
\put(8.33,.42){$\scr 2n-1$}\put(8.33,.68){$\scr 2$}
\put(8.33,1.92){$\scr n+2$}\put(8.33,2.18){$\scr n-1$}
\put(8.33,2.43){$\scr n+1$}\put(8.33,2.68){$\scr n$}

\put(4.9,3.5){\scriptsize\no de place}\put(8.3,3.5){\scriptsize\no de carte}
\put(6.95,3.55){\vector(1,0){1}}\put(7.3,3.66){$\scr h_3$}
\put(6.7,0){\line(1,0){1.5}}\put(6.7,0.25){\line(1,0){1.5}}
\put(6.7,0.5){\line(1,0){1.5}}\put(6.7,0.75){\line(1,0){1.5}}
\put(6.7,2){\line(1,0){1.5}}\put(6.7,2.25){\line(1,0){1.5}}
\put(6.7,2.5){\line(1,0){1.5}}\put(6.7,2.75){\line(1,0){1.5}}
\put(7.4,1.19){$\vdots$}

\end{picture}
\caption{Mélange de Monge, permutation $h_3$}\label{fig-Monge3}
\end{figure}
}

\newcommand{\dessinmongeshufflequatre}
{\unitlength=.1\columnwidth
\begin{figure}[h!]
\begin{picture}(10,3.8)

\put(0.7,-0.08){$\scr 1$}\put(0.7,0.43){$\scr 2$}
\put(0.22,1.97){$\scr n-1$}\put(0.7,2.43){$\scr n$}
\put(3.61,2.68){$\scr n+1$}\put(3.61,2.18){$\scr n+2$}
\put(3.61,.68){$\scr 2n-1$}\put(3.61,.18){$\scr 2n$}

\put(0.54,3.5){\scriptsize 1\ier\ paquet}
\put(2.5,3.5){\scriptsize 2\ieme\ paquet}\put(2.6,3.2){\scriptsize retourné}

\put(1,0){\line(1,0){1.5}}\put(1,0.5){\line(1,0){1.5}}
\put(1,2){\line(1,0){1.5}}\put(1,2.5){\line(1,0){1.5}}
\put(2,.25){\line(1,0){1.5}}\put(2,0.75){\line(1,0){1.5}}
\put(2,2.25){\line(1,0){1.5}}\put(2,2.75){\line(1,0){1.5}}
\put(1.7,1.07){$\vdots$}\put(2.7,1.32){$\vdots$}

\put(4.8,1.3){$\Longrightarrow$}

\put(6.4,-.08){$\scr 1$}\put(6.4,.18){$\scr 2$}
\put(6.4,.42){$\scr 3$}\put(6.4,.68){$\scr 4$}
\put(5.75,1.92){$\scr 2n-3$}\put(5.75,2.18){$\scr 2n-2$}
\put(5.75,2.42){$\scr 2n-1$}\put(6.22,2.68){$\scr 2n$}

\put(8.33,-.08){$\scr 1$}\put(8.33,.18){$\scr 2n$}
\put(8.33,.42){$\scr 2$}\put(8.33,.68){$\scr 2n-1$}
\put(8.33,1.92){$\scr n-1$}\put(8.33,2.18){$\scr n+2$}
\put(8.33,2.43){$\scr n$}\put(8.33,2.68){$\scr n+1$}

\put(4.9,3.5){\scriptsize\no de place}\put(8.3,3.5){\scriptsize\no de carte}
\put(6.95,3.55){\vector(1,0){1}}\put(7.3,3.66){$\scr h_4$}
\put(6.7,0){\line(1,0){1.5}}\put(6.7,0.25){\line(1,0){1.5}}
\put(6.7,0.5){\line(1,0){1.5}}\put(6.7,0.75){\line(1,0){1.5}}
\put(6.7,2){\line(1,0){1.5}}\put(6.7,2.25){\line(1,0){1.5}}
\put(6.7,2.5){\line(1,0){1.5}}\put(6.7,2.75){\line(1,0){1.5}}
\put(7.4,1.19){$\vdots$}

\end{picture}
\caption{Mélange de Monge, permutation $h_4$}\label{fig-Monge4}
\end{figure}
}

\newcommand{\dessinshuffleimpair}
{\unitlength=.1\columnwidth
\begin{figure}[h!]
\begin{picture}(10,3.8)

\put(0.7,0.18){$\scr 1$}\put(0.7,0.68){$\scr 2$}
\put(0.22,2.18){$\scr n-1$}\put(0.7,2.68){$\scr n$}
\put(3.61,-.08){$\scr n+1$}\put(3.61,0.42){$\scr n+2$}
\put(3.61,0.92){$\scr n+3$}\put(3.61,1.92){$\scr 2n-1$}
\put(3.61,2.42){$\scr 2n$}\put(3.61,2.92){$\scr 2n+1$}

\put(0.54,3.5){\scriptsize 1\ier\ paquet}
\put(2.54,3.5){\scriptsize 2\ieme\ paquet}
\put(1,0.25){\line(1,0){1.5}}\put(1,0.75){\line(1,0){1.5}}
\put(1,2.25){\line(1,0){1.5}}\put(1,2.75){\line(1,0){1.5}}
\put(2,0){\line(1,0){1.5}}\put(2,0.5){\line(1,0){1.5}}
\put(2,1){\line(1,0){1.5}}\put(2,2){\line(1,0){1.5}}
\put(2,2.5){\line(1,0){1.5}}\put(2,3){\line(1,0){1.5}}
\put(1.7,1.32){$\vdots$}\put(2.7,1.32){$\vdots$}

\put(6.26,-.08){$\scr 0$}\put(6.26,0.42){$\scr 1$}
\put(6.26,.92){$\scr 2$}\put(5.78,1.92){$\scr n-2$}
\put(5.78,2.42){$\scr n-1$}\put(6.26,2.92){$\scr n$}
\put(9.17,.18){$\scr n+1$}\put(9.17,0.68){$\scr n+2$}
\put(9.17,2.18){$\scr 2n-1$}\put(9.17,2.68){$\scr 2n$}

\put(6.1,3.5){\scriptsize 1\ier\ paquet}
\put(8.1,3.5){\scriptsize 2\ieme\ paquet}
\put(6.56,0){\line(1,0){1.5}}\put(6.56,0.5){\line(1,0){1.5}}
\put(6.56,1){\line(1,0){1.5}}\put(6.56,2){\line(1,0){1.5}}
\put(6.56,2.5){\line(1,0){1.5}}\put(6.56,3){\line(1,0){1.5}}
\put(7.56,0.25){\line(1,0){1.5}}\put(7.56,0.75){\line(1,0){1.5}}
\put(7.56,2.25){\line(1,0){1.5}}\put(7.56,2.75){\line(1,0){1.5}}
\put(7.26,1.32){$\vdots$}\put(8.26,1.32){$\vdots$}

\end{picture}
\caption{Jeu à un nombre impair de cartes}\label{fig-impair}
\end{figure}
}

\newcommand{\dessininshufflegene}
{\unitlength=.1\columnwidth
\begin{figure}[h!]
\begin{picture}(10,8.2)

\put(.97,1.26){$\scr 1$}
\put(.97,2.76){$\scr 2$}
\put(.97,4.26){$\scr 3$}
\put(.97,7.06){$\scr n$}

\put(2.76,.96){$\scr n+1$}
\put(2.76,2.46){$\scr n+2$}
\put(2.76,3.96){$\scr n+3$}
\put(2.87,6.76){$\scr 2n$}

\put(4.68,.66){$\scr 2n+1$}
\put(4.68,2.16){$\scr 2n+2$}
\put(4.68,3.66){$\scr 2n+3$}
\put(4.87,6.46){$\scr 3n$}

\put(8.24,.06){$\scr (k-1)n+1$}
\put(8.24,1.56){$\scr (k-1)n+2$}
\put(8.24,3.06){$\scr (k-1)n+3$}
\put(8.77,5.86){$\scr kn$}

\put(0.22,7.8){\scriptsize 1\ier\ paquet}
\put(2.28,7.8){\scriptsize 2\ieme\ paquet}
\put(4.28,7.8){\scriptsize 3\ieme\ paquet}
\put(8.18,7.8){\scriptsize $k$\ieme\ paquet}

\put(8.2,0){\line(1,0){1.5}}\put(4.3,0.6){\line(1,0){1.5}}
\put(2.3,0.9){\line(1,0){1.5}}\put(.3,1.2){\line(1,0){1.5}}
\put(8.2,1.5){\line(1,0){1.5}}\put(4.3,2.1){\line(1,0){1.5}}
\put(2.3,2.4){\line(1,0){1.5}}\put(.3,2.7){\line(1,0){1.5}}
\put(8.2,3){\line(1,0){1.5}}\put(4.3,3.6){\line(1,0){1.5}}
\put(2.3,3.9){\line(1,0){1.5}}\put(.3,4.2){\line(1,0){1.5}}
\put(8.2,5.8){\line(1,0){1.5}}\put(4.3,6.4){\line(1,0){1.5}}
\put(2.3,6.7){\line(1,0){1.5}}\put(.3,7){\line(1,0){1.5}}

\dashline{.05}(.3,0)(8.2,0)\dashline{.05}(.3,.6)(4.27,.6)
\dashline{.05}(.3,.9)(2.2,.9)\dashline{.05}(.3,1.5)(8.2,1.5)
\dashline{.05}(.3,2.1)(4.27,2.1)\dashline{.05}(.3,2.4)(2.2,2.4)
\dashline{.05}(.3,3)(8.2,3)\dashline{.05}(.3,3.6)(4.27,3.6)
\dashline{.05}(.3,3.9)(2.2,3.9)\dashline{.05}(.3,5.8)(8.2,5.8)
\dashline{.05}(.3,6.4)(4.27,6.4)\dashline{.05}(.3,6.7)(2.2,6.7)

\put(.98,5.27){$\vdots$}
\put(2.98,4.97){$\vdots$}
\put(4.98,4.67){$\vdots$}
\put(8.88,4.07){$\vdots$}
\put(6.72,6.15){$\rotatebox{-15}{\dots}$}
\put(6.72,3.35){$\rotatebox{-15}{\dots}$}
\put(6.72,1.85){$\rotatebox{-15}{\dots}$}
\put(6.72,.35){$\rotatebox{-15}{\dots}$}
\put(6.73,7.8){$\dots$}

\end{picture}
\caption{In-shuffle généralisé}\label{fig-in-shuffle-generalise}
\end{figure}
}

\newcommand{\dessininshufflepermutations}
{\unitlength=.1\columnwidth
\begin{figure}[h!]
\begin{picture}(10,6.6)
\put(1.56,-.06){$\scr 0$}\put(1.14,.44){$\scr k-2$}
\put(1.14,.69){$\scr k-1$}\put(1.56,.93){$\scr k$}
\put(1.14,1.19){$\scr k+1$}\put(.98,1.69){$\scr 2k-2$}
\put(.98,1.93){$\scr 2k-1$}\put(1.4,2.19){$\scr 2k$}
\put(.98,2.44){$\scr 2k+1$}\put(.98,2.93){$\scr 3k-2$}
\put(.98,3.19){$\scr 3k-1$}\put(1.4,3.44){$\scr 3k$}
\put(.3,4.44){$\scr (n-1)k+1$}\put(.98,4.93){$\scr nk-2$}
\put(.98,5.19){$\scr nk-1$}\put(1.4,5.44){$\scr nk$}

\put(3.53,-.06){$\scr (k-1)n+1$}\put(3.53,.44){$\scr 2n+1$}
\put(3.53,.69){$\scr n+1$}\put(3.53,.93){$\scr 1$}
\put(3.53,1.19){$\scr (k-1)n+2$}\put(3.53,1.69){$\scr 2n+2$}
\put(3.53,1.93){$\scr n+2$}\put(3.53,2.19){$\scr 2$}
\put(3.53,2.44){$\scr (k-1)n+3$}\put(3.53,2.93){$\scr 2n+3$}
\put(3.53,3.19){$\scr n+3$}\put(3.53,3.44){$\scr 3$}
\put(3.53,4.44){$\scr kn$}\put(3.53,4.93){$\scr 3n$}
\put(3.53,5.19){$\scr 2n$}\put(3.53,5.44){$\scr n$}

\put(.5,6.25){\scriptsize position}\put(3.52,6.25){\scriptsize carte}
\put(2.17,6.3){\vector(1,0){1}}\put(2.54,6.41){$\scr f$}

\put(1.9,0){\line(1,0){1.5}}\put(1.9,0.5){\line(1,0){1.5}}
\put(1.9,0.75){\line(1,0){1.5}}\put(1.9,1){\line(1,0){1.5}}
\put(1.9,1.25){\line(1,0){1.5}}\put(1.9,1.75){\line(1,0){1.5}}
\put(1.9,2){\line(1,0){1.5}}\put(1.9,2.25){\line(1,0){1.5}}
\put(1.9,2.5){\line(1,0){1.5}}\put(1.9,3){\line(1,0){1.5}}
\put(1.9,3.25){\line(1,0){1.5}}\put(1.9,3.5){\line(1,0){1.5}}
\put(1.9,4.5){\line(1,0){1.5}}\put(1.9,5){\line(1,0){1.5}}
\put(1.9,5.25){\line(1,0){1.5}}\put(1.9,5.5){\line(1,0){1.5}}
\put(2.6,3.84){$\vdots$}
\put(2.67,.125){\line(0,1){.25}}\put(2.67,1.375){\line(0,1){.25}}
\put(2.67,2.625){\line(0,1){.25}}\put(2.67,4.625){\line(0,1){.25}}

\put(6.3,-.06){$\scr 0$}\put(6.3,.19){$\scr 1$}
\put(6.3,.44){$\scr 2$}\put(5.88,.93){$\scr k-1$}
\put(6.3,1.19){$\scr k$}\put(5.88,1.44){$\scr k+1$}
\put(5.88,1.69){$\scr k+2$}\put(5.72,2.19){$\scr 2k-1$}
\put(6.14,2.44){$\scr 2k$}\put(5.72,2.69){$\scr 2k+1$}
\put(5.72,2.93){$\scr 2k+2$}\put(5.72,3.44){$\scr 3k-1$}
\put(5.46,4.44){$\scr (n-1)k$}\put(5.04,4.69){$\scr (n-1)k+1$}
\put(5.04,4.93){$\scr (n-1)k+2$}\put(5.72,5.44){$\scr nk-1$}

\put(8.23,-.06){$\scr 0$}\put(8.23,.19){$\scr n$}
\put(8.23,.44){$\scr 2n$}\put(8.23,.93){$\scr (k-1)n$}
\put(8.23,1.19){$\scr 1$}\put(8.23,1.44){$\scr n+1$}
\put(8.23,1.69){$\scr 2n+1$}\put(8.23,2.19){$\scr (k-1)n+1$}
\put(8.23,2.44){$\scr 2$}\put(8.23,2.69){$\scr n+2$}
\put(8.23,2.93){$\scr 2n+2$}\put(8.23,3.44){$\scr (k-1)n+2$}
\put(8.23,4.44){$\scr n-1$}\put(8.23,4.69){$\scr 2n-1$}
\put(8.23,4.93){$\scr 3n-1$}\put(8.23,5.44){$\scr kn-1$}

\put(5.2,6.25){\scriptsize position}\put(8.2,6.25){\scriptsize carte}
\put(6.85,6.3){\vector(1,0){1}}\put(7.22,6.42){$\scr \tg$}

\put(6.6,0){\line(1,0){1.5}}\put(6.6,0.25){\line(1,0){1.5}}
\put(6.6,0.5){\line(1,0){1.5}}\put(6.6,1){\line(1,0){1.5}}
\put(6.6,1.25){\line(1,0){1.5}}\put(6.6,1.5){\line(1,0){1.5}}
\put(6.6,1.75){\line(1,0){1.5}}\put(6.6,2.25){\line(1,0){1.5}}
\put(6.6,2.5){\line(1,0){1.5}}\put(6.6,2.75){\line(1,0){1.5}}
\put(6.6,3){\line(1,0){1.5}}\put(6.6,3.5){\line(1,0){1.5}}
\put(6.6,4.5){\line(1,0){1.5}}\put(6.6,4.75){\line(1,0){1.5}}
\put(6.6,5){\line(1,0){1.5}}\put(6.6,5.5){\line(1,0){1.5}}
\put(7.3,3.84){$\vdots$}
\put(7.37,.625){\line(0,1){.25}}\put(7.37,1.875){\line(0,1){.25}}
\put(7.37,3.125){\line(0,1){.25}}\put(7.37,5.125){\line(0,1){.25}}

\end{picture}
\caption{$\!\!$Mélanges généralisés, $\!$permutations $\!f\!$ et $\!\tg$}
\label{fig-melanges-generalises}
\end{figure}
}

\newcommand{\dessinoutshufflegene}
{\unitlength=.1\columnwidth
\begin{figure}[h!]
\begin{picture}(10,7.9)

\put(.97,.06){$\scr 0$}
\put(.97,1.56){$\scr 1$}
\put(.97,3.06){$\scr 2$}
\put(.76,5.56){$\scr n-1$}

\put(2.97,.36){$\scr n$}
\put(2.76,1.86){$\scr n+1$}
\put(2.76,3.36){$\scr n+2$}
\put(2.68,5.86){$\scr 2n-1$}

\put(4.87,.66){$\scr 2n$}
\put(4.68,2.16){$\scr 2n+1$}
\put(4.68,3.66){$\scr 2n+2$}
\put(4.68,6.16){$\scr 3n-1$}

\put(8.44,1.26){$\scr (k-1)n$}
\put(8.24,2.76){$\scr (k-1)n+1$}
\put(8.24,4.26){$\scr (k-1)n+2$}
\put(8.58,6.76){$\scr kn-1$}

\put(0.22,7.5){\scriptsize 1\ier\ paquet}
\put(2.28,7.5){\scriptsize 2\ieme\ paquet}
\put(4.28,7.5){\scriptsize 3\ieme\ paquet}
\put(8.18,7.5){\scriptsize $k$\ieme\ paquet}

\put(8.2,1.2){\line(1,0){1.5}}\put(4.3,0.6){\line(1,0){1.5}}
\put(2.3,0.3){\line(1,0){1.5}}\put(.3,0){\line(1,0){1.5}}
\put(8.2,2.7){\line(1,0){1.5}}\put(4.3,2.1){\line(1,0){1.5}}
\put(2.3,1.8){\line(1,0){1.5}}\put(.3,1.5){\line(1,0){1.5}}
\put(8.2,4.2){\line(1,0){1.5}}\put(4.3,3.6){\line(1,0){1.5}}
\put(2.3,3.3){\line(1,0){1.5}}\put(.3,3){\line(1,0){1.5}}
\put(8.2,6.7){\line(1,0){1.5}}\put(4.3,6.1){\line(1,0){1.5}}
\put(2.3,5.8){\line(1,0){1.5}}\put(.3,5.5){\line(1,0){1.5}}

\dashline{.05}(.3,.3)(2.2,.3)\dashline{.05}(.3,.6)(4.27,.6)
\dashline{.05}(.3,1.2)(8.2,1.2)\dashline{.05}(.3,1.8)(2.2,1.8)
\dashline{.05}(.3,2.1)(4.27,2.1)\dashline{.05}(.3,2.7)(8.2,2.7)
\dashline{.05}(.3,3.3)(2.2,3.3)\dashline{.05}(.3,3.6)(4.27,3.6)
\dashline{.05}(.3,4.2)(8.2,4.2)\dashline{.05}(.3,5.8)(2.2,5.8)
\dashline{.05}(.3,6.1)(4.27,6.1)\dashline{.05}(.3,6.7)(8.2,6.7)

\put(.98,4.37){$\vdots$}
\put(2.98,4.67){$\vdots$}
\put(4.98,4.97){$\vdots$}
\put(8.88,5.57){$\vdots$}
\put(6.72,6.31){$\rotatebox{15}{\dots}$}
\put(6.72,3.81){$\rotatebox{15}{\dots}$}
\put(6.72,2.31){$\rotatebox{15}{\dots}$}
\put(6.72,.81){$\rotatebox{15}{\dots}$}
\put(6.73,7.5){$\dots$}

\end{picture}
\caption{Out-shuffle généralisé}\label{fig-out-shuffle-generalise}
\end{figure}
}

\begin{document}
\maketitle

\section{\textsf{Description du problème}}

Le problème suivant m'a été soumis par l'un de mes
élèves\footnote{Anthony Tschirhard, INSA de Lyon, 51\ieme\ promotion} qui a
une grande pratique de la magie. On dispose d'une pile de $2n$ jetons dont les
$n$~du bas sont verts et les $n$ du haut sont rouges. Il divise la pile en
deux parties : l'une contenant les $n$~jetons verts, l'autre contenant les $n$
jetons rouges. Avec une dextérité remarquable, de sa seule main droite il
prend simultanément les deux piles de jetons, puis procède à un mélange
de ces deux dernières pour former une unique pile (de $2n$ jetons) alternant
parfaitement jeton rouge et jeton vert. Cette nouvelle pile est divisée
à son tour exactement en son milieu donnant deux tas de $n$ jetons et notre
magicien procède à un nouveau mélange intercalant alternativement un jeton de
chaque tas pour former une nouvelle pile de $2n$ jetons qu'il recoupe de
nouveau en deux piles de $n$ jetons et ainsi de suite.

Après plusieurs expériences (pour $n\le 8$), notre magicien s'aperçut qu'il
retombait au bout d'un certain nombre de tels mélanges sur la pile initiale :
les $n$ jetons du bas sont verts et les $n$ du haut sont rouges
(voir Fig.~\ref{fig-jetons1}--\ref{fig-jetons6} pour $n=4$).
Plus frappant encore, en numérotant et ordonnant les jetons, il semblerait
que la première fois que l'on retrouve une pile de $n$ jetons verts consécutifs
et de $n$~rouges de bas en haut, les jetons soient en fait exactement
dans l'ordre initial. Les questions qu'il m'a soumises étaient
alors les suivantes : étant donnée une pile de $2n$ jetons numérotés
$1,2,\dots,2n-1,2n$ ($n\in\N^*$) et superposés de bas en haut dans l'ordre
croissant,
\begin{enumerate}
\item
le processus de mélange décrit ci-dessus conduit-il nécessairement au classement
initial en un temps fini ?
\item
si oui, peut-on exprimer en fonction de $n$ le nombre minimum
de mélanges nécessaires pour retrouver le classement initial ?
\item
si les $n$ premiers jetons sont verts et les $n$~derniers rouges,
peut-on arriver au regroupement initial des jetons (la pile du bas constituée de
jetons verts, celle du haut de rouges) éventuellement dans le désordre quant
à la numérotation ?
\end{enumerate}

Il s'agit d'un problème classique bien connu du monde de la magie des cartes,
la manipulation précédente étant souvent faite avec des cartes à jouer.
Le jeu habituel de $32$ cartes est priviligié car le nombre $32$ se décompose
en $2^5$, ce qui conduit à des propriétés remarquables.
Dans la littérature, ce type de problème est abordé sous le vocable anglophone
de <<~chip-shuffle~>> pour les jetons, de <<~riffle-shuffle~>> pour les cartes
et plus précisément de <<~in-shuffle~>> et de <<~out-shuffle~>>
selon le placement des cartes initiales de chaque paquet lors d'un mélange.
Il a été considéré entre autres par le célèbre informaticien et magicien Elmsley
en 1957 (\cite{elmsley}). On trouvera également une large description historique
du problème dans la section intitulée <<~some history of perfect shuffle~>>
de~\cite{diaconis-graham-kantor}.

Les in-shuffles et out-shuffles sont des mélanges très proches.
Disposant d'un jeu de cartes que l'on coupe au milieu donnant deux paquets
de cartes (un premier et un deuxième), l'in-shuffle du jeu initial consiste à
démarrer le processus en plaçant la première carte du bas du premier paquet
sur la première carte du bas du deuxième paquet. \`A l'opposé, l'out-shuffle
consiste à démarrer en plaçant la première carte du bas du premier paquet
sous la première carte du bas du deuxième paquet.
Notons que dans le cas de l'out-shuffle, les première et dernière cartes
restent immobiles tout au long de l'expérience.
On pourra consulter les sites internet \cite{wiki,wolfram} pour un recensement
de ces divers mélanges ainsi que d'autres.

\dessinjetonsun

\dessinjetonsdeux

\dessinjetonstrois

\dessinjetonsquatre

\dessinjetonscinq

\dessinjetonssix

Elmsley se posait la question plus complexe : est-il possible de déplacer
la carte du dessus du paquet à une position donnée à l'avance à l'aide d'une
succession de mélanges de nature in-shuffle ou out-shuffle (\cite{elmsley})~?
La réponse est positive et il obtint un procédé pour accomplir une telle
manipulation. Inversement, est-il possible de faire apparaître une carte
donnée dans le jeu sur le haut du paquet de cartes, voire
à une place quelconque selon un précédé similaire ?
Récemment, les mathématiciens Diaconis et Graham (le premier auteur est
également magicien) ont répondu positivement à la question générale du
déplacement d'une carte donnée vers une position donnée.
Ils ont proposé un algorithme indiquant le chemin à suivre (succession de
in- et out-shuffles, \cite{diaconis-graham}).

Reprenons les questions posées par mon élève.
La réponse à la première question est affirmative et la raison en est très simple.
Sommairement, on travaille dans un groupe fini de permutations, on revient
donc à la position initiale au bout d'un nombre fini d'expériences, il s'agit
d'un processus périodique.
La réponse à la deuxième question est également affirmative ; c'est un calcul
de période et une formule implicite est disponible, voir \eg\ les livres de
Conway \& Guy, \cite[p. 163--165]{conway}, de Herstein \& Kaplansky,
\cite[Chap. 3.4]{kaplan}, ou d'Uspensky \& Heaslet, \cite[p. 244--245]{uspensky}.
Néanmoins, une formule explicite ne semble pas accessible excepté pour
les cas particuliers où $n$ est une puissance de $2$ à $\pm1$ près.
La troisième question semble rester à ma connaissance ouverte.

Dans cet article, je détaille le processus de l'in-shuffle -- celui de
l'out-shuffle s'en déduisant aisément -- à l'aide de permutations de l'ensemble
$\{0,1,2,\dots,2n-1\}$ ou $\{1,2,3,\dots,2n\}$. L'un des deux ensembles
sera d'une utilisation plus commode que l'autre selon le cas étudié.
Les permutations associées à $\{0,1,2,\dots,2n-1\}$ sont particulièrement
bien adaptées au calcul explicite des itérations successives
(correspondant à la succession de mélanges parfaits) par le truchement
des écritures binaires.
Aussi, je détaille toutes les permutations relatives aux deux numérotations
en signalant laquelle permet de formuler le plus simplement possible la période
recherchée. Je porte un intérêt particulier aux cas spécifiques mentionnés
plus haut, à savoir lorsque $n$ est une puissance de $2$ à $\pm1$ près.
J'examine également d'autres exemples de mélanges parfaits : les mélanges
de Monge qui sont des in/out-shuffles déformés par une symétrie (voir \eg\
\cite{conway,wiki,wolfram}).
Toute cette analyse permet d'accéder au calcul
des périodes de chacun des mélanges décrits ici.

Par souci d'homogénéité des notations, on notera $f,g,h,\dots$ les permutations
relatives à l'énumération $1,2,3,\dots,2n$ et $\tf,\tg,\th,\dots$ celles
associées à $0,1,2,\dots,2n-1$. On passe par exemple de la permutation $f$ à
la permutation $\tf$ selon
la relation $\tf(i)=f(i+1)-1$ pour tout \mbox{$i\in\{0,1,2,\dots,2n-1\}$}.

Il est remarquable de constater que, bien au-delà de son aspect ludique, ce
problème suscite des questions de nature algébrique avancées (théorie des
groupes, \cite{diaconis-graham-kantor,golomb}). L'étude approfondie de certaines
permutations mises en jeu a été entreprise indépendamment de ce contexte par
Lévy (\cite{levy1,levy2,levy3}).
Ce problème connaît également des applications notamment en informatique
(calcul parallèle, réseaux, \cite{chen,schwartz,stone}). Mentionnons enfin
l'existence d'autres types de mélanges d'une importance notoire qui ont retenu
l'attention des mathématiciens : les mélanges aléatoires (voir~\cite{biane}
pour plus d'informations concernant ce type de mélange).

\vspace{\baselineskip}
\centerline{\textbf{\textsf{ Plan de l'article}}}
\begin{itemize}
\item
Dans la section~\ref{section-model}, nous présentons la modélisation du
problème en introduisant toutes les permutations relatives aux différents
mélanges de cartes. Certaines permutations redondantes ne seront pas
nécessairement utiles, mais nous avons choisi de les écrire systématiquement
afin de garder une homogénéité de notations ainsi qu'une logique de présentation.
Selon les cas étudiés, certaines d'entre elles serviront à calculer la période
du mélange en question, alors que d'autres permettront le calcul
explicite des itérations successives.
\item
Dans la section~\ref{section-solut}, nous formulons de manière implicite
les réponses relatives au calcul de périodes des mélanges étudiés.
\item
Dans les sections~\ref{section-caspart1} et~\ref{section-caspart2},
nous calculons explicitement, en faisant appel au calcul binaire, toutes les
itérations successives de certaines permutations significatives, dans les cas
particuliers où $n$ est une puissance de $2$ à $\pm1$ près,
ce qui permet de fournir une réponse explicite aux questions posées.
\item
Dans la section~\ref{section-impair}, nous examinons brièvement le cas
d'un jeu contenant un nombre impair de cartes.
\item
Dans la section~\ref{section-deplacement}, nous décrivons un procédé simple
pour déplacer une carte donnée vers une position prédéfinie dans le cas où
$n$ est une puissance de $2$.
\item
Nous terminons dans la section~\ref{section-gene} par une tentative de
généralisation à un paquet de cartes divisible par un nombre supérieur à deux.
\end{itemize}

\section{\textsf{Modélisation}}\label{section-model}

Le jeu de cartes numérotées de bas en haut dans l'ordre $1,2,$ $3,\dots,2n$ est coupé en son
milieu et donne les deux paquets de cartes numérotées
$1,2,\dots,n$ pour le premier et $n+1,n+2,\dots,2n$ pour le deuxième.
On constitue un nouveau jeu de $2n$ cartes en prenant à tour de rôle
une carte de chacun des paquets de $n$ cartes.

\subsection{\textsf{In-shuffle}}

Si l'on démarre le mélange par le deuxième paquet (cas d'un in-shuffle), on
obtient dans l'ordre la suite de cartes
\nos\ $n+1,1,n+2,2,n+3,3,\dots,\lb 2n-1,n-1,2n,n$ (voir Fig.~\ref{fig-in-shuffle}).

\dessininshuffle

Cela se traduit par une correspondance entre numéros d'ordre avant et après
mélange : à l'issue du mélange, la première carte porte le \no $n+1$, la
deuxième le \no $1$, la troisième le \no $n+2$, la quatrième le \no $2$, etc.
De manière générale, la $i$\ieme\ carte porte le \no $i/2$ lorsque $i$ est pair
et le \no $(i+1)/2+n$ lorsque $i$ est impair. Cette correspondance définit une
permutation $f$ des entiers $1,2,3,\dots,2n$ donnant le numéro de la carte se
situant à la $i\ieme$\ position à l'issue du mélange :
$$
f(i)=\left\{\begin{array}{ll}
\dis\frac{i}{2} & \mbox{si $i$ est pair,}
\\[1.5ex]
\dis\frac{i+1}{2}+n & \mbox{si $i$ est impair,}
\end{array}\right.
$$
que l'on peut encore écrire $f(i)=[\frac{i+1}{2}]+\e(i)n$ où $[\,]$ est la
fonction partie-entière et $\e(i)=0$ si $i$ est pair, $\e(i)=1$ si $i$ est impair.

La permutation réciproque est un peu plus simple à écrire :
$$
f^{-1}(j)=\left\{\begin{array}{ll}
2j & \mbox{si $1\le j\le n$,}
\\
2j-2n-1 & \mbox{si $n+1\le j\le 2n$.}
\end{array}\right.
$$
Rappelant la définition d'une congruence arithmétique,
<<~$\modtext{a}{b}{n}$~>> signifie <<~$a-b$ est divisible par $n$~>>,
on a en particulier la congruence remarquable
$$
\modulo{f^{-1}(j)}{2j}{(2n+1)}
$$
qui sera utile pour calculer la période de $f$.
Cette réciproque indique qu'une carte portant un numéro $j$ entre $1$ et
$n$ occupe à l'issue du mélange la place \no $2j$, et qu'une carte
portant un numéro $j$ entre $n+1$ et $2n$ occupe la place \no  $2j-2n-1$.
Par exemple, la carte \no $1$ occupe la place \no $2$, la carte \no $2$ occupe
la place \no $4,\dots$, puis la carte \no $n$ occupe la place \no $2n$,
la carte \no $n+1$ occupe la place \no $1$, la carte \no $n+2$ occupe
la place \no $3$, etc.

Il est utile de transcrire cette modélisation en translatant simplement la
numérotation des cartes d'une unité. Cela fournit la permutation $\tf$ des
entiers $0,1,2,\dots,2n-1$ définie par \mbox{$\tf(i)=f(i+1)-1$,} soit :
$$
\tf(i)=\left\{\begin{array}{ll}
\dis\frac{i}{2}+n & \mbox{si $i$ est pair,}
\\[1.5ex]
\dis\frac{i-1}{2} & \mbox{si $i$ est impair,}
\end{array}\right.
$$
et la permutation réciproque s'écrit :
$$
\tf^{-1}(j)=\left\{\begin{array}{ll}
2j+1 & \mbox{si $0\le j\le n-1$,}
\\
2j-2n & \mbox{si $n\le j\le 2n-1$.}
\end{array}\right.
$$

Les permutations $f$ et $\tf$ nous permettront de calculer explicitement la
période de l'in-shuffle dans les cas particuliers $n=2^{p-1}$ et $n=2^{p-1}-1$.

\subsection{\textsf{Out-shuffle}}

Si on démarre le mélange à présent par le premier paquet (cas d'un out-shuffle),
on obtient dans l'ordre les cartes \nos\
$1,n+1,2,n+2,3,\lb n+3,\dots,n-1,2n-1,n,2n$ (voir Fig.~\ref{fig-out-shuffle}).

\dessinoutshuffle

Ce processus définit alors la permutation $g$ des entiers $1,2,\dots,2n$ suivante :
$$
g(i)=\left\{\begin{array}{ll}
\dis\frac{i}{2}+n & \mbox{si $i$ est pair,}
\\[1.5ex]
\dis\frac{i+1}{2} & \mbox{si $i$ est impair,}
\end{array}\right.
$$
de réciproque
$$
g^{-1}(j)=\left\{\begin{array}{ll}
2j-1 & \mbox{si $1\le j\le n$,}
\\
2j-2n & \mbox{si $n+1\le j\le 2n$.}
\end{array}\right.
$$
On observe que $g(1)=1$ et $g(2n)=2n$, \ie\ les cartes \nos\ $1$ et $2n$ restent
immobiles tout au cours de la manipulation.

Retranscrivons cette modélisation en translatant la numérotation des cartes
d'une unité. Cela fournit la permutation $\tg$ des entiers
$0,1,2,\dots,2n-1$ définie par $\tg(i)=g(i+1)-1$, soit :
$$
\tg(i)=\left\{\begin{array}{ll}
\dis\frac{i}{2} & \mbox{si $i$ est pair,}
\\[1.5ex]
\dis\frac{i-1}{2}+n & \mbox{si $i$ est impair,}
\end{array}\right.
$$
et la permutation réciproque s'écrit :
$$
\tg^{-1}(j)=\left\{\begin{array}{ll}
2j & \mbox{si $0\le j\le n-1$,}
\\
2j+1-2n & \mbox{si $n\le j\le 2n-1$.}
\end{array}\right.
$$
Notons en particulier la congruence
$$
\modulo{\tg^{-1}(j)}{2j}{(2n-1)}.
$$
La restriction de $\tg$ à l'ensemble $\{1,\dots,2n-2\}$ est une permutation
qui correspond à un in-shuffle de $2n-2$ cartes.
Cette observation indique qu'un out-shuffle de $2n$ cartes est identique à
un in-shuffle de ce jeu auquel on a retiré les première et dernière cartes,
donnant ainsi un jeu de $2n-2$ cartes.

\subsection{\textsf{Mélange de Monge : première version}}

Le mélange de Monge consiste à faire passer les cartes d'un jeu complet
d'une main à l'autre en intercalant alternativement les cartes l'une au-dessus
et au-dessous de l'autre.

\subsubsection{\textsf{Première façon de démarrer le processus}}\label{monge-1}

On numérote les cartes de $1$ à $2n$. On commence par prendre la carte \no $1$
qui va démarrer le nouveau paquet, puis la carte \no $2$ que l'on place
au-dessous du nouveau paquet (donc au-dessous de la carte \no $1$), puis la
carte \no $3$ que l'on place au-dessus du nouveau paquet (donc au-dessus de la
\no $1$), puis la carte \no $4$ que l'on place au-dessous du nouveau paquet
(donc au-dessous de la \no $2$), puis la carte \no $5$ que l'on place au-dessus
du nouveau paquet (donc au-dessus de la \no $3$), et ainsi de suite
(voir Fig.~\ref{fig-Monge1}).
On obtient alors dans l'ordre, de bas en haut, les cartes \nos\
$2n,2n-2, 2n-4,\dots,4,2,1,3,\dots,2n-3,2n-1$.
On rencontre donc les cartes de numéros pairs dans l'ordre décroissant puis
celles de numéros impairs dans l'ordre croissant.

\dessinmongeshuffleun

Ce mélange correspond à la
permutation $h_1^{-1}$ des entiers $1,2,\dots,2n$ où $h_1$ est donnée par
$$
h_1(i)=\left\{\begin{array}{ll}
\dis\frac{i+1}{2}+n & \mbox{si $i$ est impair,}
\\[1ex]
\dis n+1-\frac{i}{2} & \mbox{si $i$ est pair.}
\end{array}\right.
$$
La permutation $h_1^{-1}$ s'écrit
$$
h_1^{-1}(j)=\left\{\begin{array}{ll}
2n+2-2j & \mbox{si $1\le j\le n$,}
\\
2j-2n-1 & \mbox{si $n+1\le j\le 2n$.}
\end{array}\right.
$$
La permutation analogue $\th_1$ associée à la numérotation $0,1,2,\dots,2n-1$
est donnée par
$$
\th_1(i)=\left\{\begin{array}{ll}
\dis\frac{i}{2}+n & \mbox{si $i$ est pair,}
\\[1ex]
\dis n-\frac{i+1}{2} & \mbox{si $i$ est impair,}
\end{array}\right.
$$
et sa réciproque par
$$
\th_1^{-1}(j)=\left\{\begin{array}{ll}
2n-1-2j & \mbox{si $0\le j\le n-1$,}
\\
2j-2n & \mbox{si $n\le j\le 2n-1$.}
\end{array}\right.
$$

\subsubsection{\textsf{Deuxième façon de démarrer le processus}}\label{monge-2}

On commence par prendre la carte \no $1$ qui va démarrer le nouveau paquet,
puis la carte \no $2$ que l'on place à présent au-dessus du nouveau paquet
(donc au-dessus de la carte \no $1$), puis la carte \no $3$ que l'on place
au-dessous de la \no $1$, puis la carte \no $4$ que l'on place au-dessus de la
\no $2$, puis la carte \no $5$ que l'on place au-dessous de la \no $3$,
et ainsi de suite (voir Fig.~\ref{fig-Monge2}).
On obtient cette fois dans l'ordre, de bas en haut, les cartes
\nos\ $2n-1,\lb 2n-3,\dots,3,1,2,4\dots,2n-2,2n$,
c'est-à-dire les cartes de numéros impairs dans l'ordre décroissant puis
celles de numéros pairs dans l'ordre croissant.

\dessinmongeshuffledeux

Remarquons que ce mélange se déduit du précédent simplement en retournant
le jeu, ou encore en effectuant la symétrie $j\longmapsto 2n+1-j$.\lb
En notant $h_2^{-1}$ la permutation des entiers $1,2,\dots,2n$ correspondant à
ce mélange, on a précisément la relation $h_2^{-1}(j)=h_1^{-1}(2n+1-j)$ pour tout
$j\in\{1, 2,\dots,2n\}$, qui est équivalente à \lb $h_2(i)=2n+1-h_1(i)$
pour tout $i\in\{1,2,\dots,2n\}$. Donc
$$
h_2(i)=\left\{\begin{array}{ll}
\dis\frac{i}{2}+n & \mbox{si $i$ est pair,}
\\[1ex]
\dis n-\frac{i-1}{2} & \mbox{si $i$ est impair,}
\end{array}\right.
$$
et
$$
h_2^{-1}(j)=\left\{\begin{array}{ll}
2n+1-2j & \mbox{si $1\le j\le n$,}
\\
2j-2n & \mbox{si $n+1\le j\le 2n$.}
\end{array}\right.
$$
La permutation $\th_2$ relative à la numérotation $0,1,2,\dots,2n-1$ s'écrit :
$$
\th_2(i)=\left\{\begin{array}{ll}
\dis\frac{i-1}{2}+n & \mbox{si $i$ est impair,}
\\[1ex]
\dis n-1-\frac{i}{2} & \mbox{si $i$ est pair,}
\end{array}\right.
$$
et sa réciproque :
$$
\th_2^{-1}(j)=\left\{\begin{array}{ll}
2n-2-2j & \mbox{si $0\le j\le n-1$,}
\\
2j-2n+1 & \mbox{si $n\le j\le 2n-1$.}
\end{array}\right.
$$

\subsection{\textsf{Mélange de Monge : deuxième version}}

Une autre version du mélange de Monge consiste à couper le jeu de cartes
initialement numérotées $1,2,\dots,2n$ en deux paquets de cartes numérotées
$1,2,\dots,n$ et $n+1,n+2,\dots,2n$, à retourner le deuxième paquet qui devient
ordonné selon $2n,2n-1,\dots,n+2,n+1$, puis de mélanger le premier paquet
et le deuxième ainsi retourné via un in-shuffle ou un out-shuffle. Cela donne
un mélange <<~à l'espagnole~>> : on dispose les deux paquets sous forme
d'éventails, on retourne le deuxième éventail (Olé !) et on intercale
parfaitement les deux éventails.

\subsubsection{\textsf{In-shuffle de Monge}}\label{sous-section}

On intercale les deux paquets de cartes numérotées $1,2,\dots,n$ et
$2n,2n-1,\dots,n+2,n+1$ en commençant par le deuxième (cas d'un in-shuffle,
voir Fig.~\ref{fig-Monge3}).
Cela donne la succession de cartes \nos\
$2n,1,2n-1,2,\dots, n+2,n-1,n+1,n.$

\dessinmongeshuffletrois

Cette séquence est représentée par la
permutation $h_3$ des entiers $1,2,\dots,2n$ suivante :
$$
h_3(i)=\left\{\begin{array}{ll}
\dis\frac{i}{2} & \mbox{si $i$ est pair,}
\\
\dis 2n-\frac{i-1}{2} & \mbox{si $i$ est impair,}
\end{array}\right.
$$
de réciproque
$$
h_3^{-1}(j)=\left\{\begin{array}{ll}
2j & \mbox{si $1\le j\le n$,}
\\
4n+1-2j & \mbox{si $n+1\le j\le 2n$.}
\end{array}\right.
$$
On a en particulier la congruence importante
$$
\modulo{h_3^{-1}(j)}{\pm 2j}{(4n+1)}.
$$
Les permutations analogues des entiers $0,1,2,\dots,\linebreak 2n-1$ sont
alors les suivantes :
$$
\th_3(i)=\left\{\begin{array}{ll}
\dis\frac{i-1}{2} & \mbox{si $i$ est impair,}
\\
\dis 2n-1-\frac{i}{2} & \mbox{si $i$ est pair,}
\end{array}\right.
$$
et
$$
\th_3^{-1}(j)=\left\{\begin{array}{ll}
2j+1 & \mbox{si $0\le j\le n-1$,}
\\
4n-2-2j & \mbox{si $n\le j\le 2n-1$.}
\end{array}\right.
$$
Ce mélange est relié à celui décrit dans~\ref{monge-1} selon la relation
$2n+1-h_3(2n+1-i)=h_1(i)$ valable pour tout $i\in\{1,2,\dots,2n\}$,
soit encore en notant $s$ la symétrie $i\longmapsto 2n+1-i$,
$$
h_3=s\circ h_1\circ s=s\circ h_1\circ s^{-1}.
$$
Les permutations $h_1$ et $h_3$ sont donc conjuguées dans le groupe des
permutations $(\cS_{2n},\circ)$.

\subsubsection{\textsf{Out-shuffle de Monge}}

On intercale à présent les deux paquets de cartes numérotées
$1,2,\dots,n$ et $2n,2n-1,\dots,\lb n+2,n+1$ en commençant par
le premier (cas d'un out-shuffle, voir Fig.~\ref{fig-Monge4}).
Cela conduit à la séquence de cartes \nos\ $1,2n,2,2n-1,\dots,n-1,n+2,n,n+1.$

\dessinmongeshufflequatre

D'où les permutations des entiers $1,2,\dots,2n$ suivantes :
$$
h_4(i)=\left\{\begin{array}{ll}
\dis\frac{i+1}{2} & \mbox{si $i$ est impair,}
\\
\dis 2n+1-\frac{i}{2} & \mbox{si $i$ est pair,}
\end{array}\right.
$$
et
$$
h_4^{-1}(j)=\left\{\begin{array}{ll}
2j-1 & \mbox{si $1\le j\le n$,}
\\
4n+2-2j & \mbox{si $n+1\le j\le 2n$.}
\end{array}\right.
$$
On a $h_4(1)=1$, ce qui signifie que la carte \no 1 reste immobile dans cette
manipulation. Les permutations des entiers $0,1,2,\dots,2n-1$ associées sont
les suivantes :
$$
\th_4(i)=\left\{\begin{array}{ll}
\dis\frac{i}{2} & \mbox{si $i$ est pair,}
\\
\dis 2n-\frac{i+1}{2} & \mbox{si $i$ est impair,}
\end{array}\right.
$$
et
$$
\th_4^{-1}(j)=\left\{\begin{array}{ll}
2j & \mbox{si $0\le j\le n-1$,}
\\
4n-1-2j & \mbox{si $n\le j\le 2n-1$.}
\end{array}\right.
$$
On a en particulier la congruence notoire
$$
\modulo{\th_4^{-1}(j)}{\pm 2j}{(4n-1)}.
$$
\`A l'instar du mélange de la section~\ref{sous-section},
ce mélange est relié à celui décrit dans~\ref{monge-2} selon
$$
h_4=s\circ h_2\circ s^{-1}.
$$
Les permutations $h_2$ et $h_4$ sont donc conjuguées.

Pour résumer et avoir une vision concise des divers mélanges considérés, nous
avons défini les permutations $f$ (in-shuffle), $g$ (out-shuffle) et
$h_1,h_2,h_3,h_4$ (mélanges de Monge) ainsi que leurs translatées
$\tf,\tg,\th_1,\th_2,\th_3,\th_4$ et leurs réciproques.

\section{\textsf{Une formulation de la solution du problème}}\label{section-solut}

Les mélanges répétés correspondent mathématiquement aux itérations successives
des permutations décrites précédemment. Travaillons par exemple avec la
permutation de l'in-shuffle $f$. L'issue des $k$ premiers mélanges est
représentée par la permutation
$$
f^k=\underset{k \mbox{ \scriptsize fois}}{\underbrace{f\circ\dots\circ f}}.
$$
Plus précisément, la quantité $f^k(i)$ désigne le numéro de la carte se
trouvant à la position \no $i$.
Ainsi l'intégralité de l'expérience (supposée illimitée...) est modélisée par
l'ensemble des itérations successives de $f$ suivant :
$\{id,f,f^2,f^3,\dots\}\lb =\{f^k,k\in\N\}.$ C'est un sous-ensemble de l'ensemble
fini $\cS_{2n}$ des permutations de $\{1,2,\dots,2n\}$, il est donc fini
lui-même et il y a au moins deux entiers distincts $k$ et $l$ tels que
$f^k\neq f^l$. Comme $f$~est bijective, on a en supposant par exemple que
$k<l$ et en posant alors $r=l-k$, $f^r=id$, ou encore $(f^{-1})^r=id$.
Cela signifie qu'au bout de $r$~mélanges parfaits, on retombe nécessairement
sur l'ordre initial des cartes. Ceci répond affirmativement à la première
question posée dans l'introduction.

L'ensemble $\{f^k,k\in\N\}$, qui est \textit{a posteriori} identique à
$\{f^k,k\in\Z\}$, s'écrit en extension selon $\{id,f,f^2,\dots,f^{r-1}\}$.
Le nombre $r$ est une période de l'application $k\in\Z\longmapsto f^k\in\cS_{2n}$.
Remarquons que la permutation $\tf$ qui est définie par $\tf(i)=f(i+1)-1$ vérifie
$\tf^k(i)=f^k(i+1)-1$.
Ainsi, l'équation d'inconnue $r$, $f^r=id$ est équivalente à $\tf^r=id$.
Cela signifie que les permutations $f$ et $\tf$ (et aussi $f^{-1}$ et $\tf^{-1}$)
ont même période, ce qui est bien sûr naturel puisque le retour à l'ordre initial
ne dépend pas du choix de la numérotation.
Nous pourrons travailler indifféremment avec $f$ ou $\tf$ selon le cas.
Le calcul de cette période -- le plus petit $r\ge 1$ tel que $f^r=id$ --
est précisément l'objet de la deuxième question posée dans l'introduction.
Une méthode de calcul est proposée dans les théorèmes~\ref{th-in-shuffle},
\ref{th-out-shuffle}, \ref{th-monge-shuffle1} et~\ref{th-monge-shuffle2}
ci-dessous, on peut la trouver \eg\ dans les divers livres \cite{conway,kaplan,uspensky}
ainsi que dans les articles \cite{elmsley,morris1}.

La dernière question posée dans l'introduction concerne la possibilité de
retrouver, partant d'une pile de $2n$ jetons contenant de bas en haut
$n$ jetons verts et $n$ jetons rouges, une pile visuellement identique, mais
éventuellement correspondant à un ordre de jetons (s'ils étaient discernables)
différent. Cela revient à déterminer pour la permutation $f$
le plus petit entier $r'\ge 1$ tel que
$$
f^{r'}(\{1,2,\dots,n\})=\{1,2,\dots,n\}
\vspace{-\baselineskip}
$$
et
$$
f^{r'}(\{n+1,n+2,\dots,2n\})=\{n+1,n+2,\dots,2n\}.
$$
Une des deux égalités est redondante puisque $f$ est bijective.
Il est clair que $r'\le r$. En fait, l'entier $r'$ divise la période $r$.
En effet, en effectuant la division euclidienne de $r$ par $r'$, on~a
$r=ar'+b$ pour deux entiers $a$ et $b$ tels que $0\le b\le r'-1$.
Alors, par définition de $r$, $f^r=f^{(a+1)r'-(r'-b)}=id$, donc
$f^{r'-b}=(f^{r'})^{a+1}$ et $f^{r'-b}$ vérifie
$f^{r'-b}(\{1,2,\dots,n\})=\{1,2,\dots,n\}$.
Par définition de $r'$, comme $1\le r'-b\le r'$, on a nécessairement $b=0$
et donc $r'$ divise $r$. La recherche de $r'$ semble délicate, nous
ne l'aborderons pas dans ce travail.

\subsection{\textsf{Cas de l'in-shuffle}}

\bth{\label{th-in-shuffle}
La période de l'in-shuffle de $2n$ cartes est le plus petit entier $r\ge 1$
vérifiant $\modtext{2^r}{1}{(2n+1)}$. En d'autres termes, $r$ est l'ordre de
$2$ modulo $(2n+1)$.
}
\dem
Rappelons la congruence pour la permutation $f$ associée à l'in-shuffle :
$$
\modulo{f^{-1}(j)}{2j}{(2n+1)}.
$$
On a alors pour tout $k\in\N$,
$$
\modulo{f^{-k}(j)}{2^kj}{(2n+1)}.
$$
La période de $f$ est donc le plus petit entier $r\ge 1$ vérifiant
$\modtext{2^r}{1}{(2n+1)}$. Notons qu'un tel entier existe effectivement. En
effet, en adaptant le raisonnement précédent et en rappelant la notation
$\modsimple{a}{n}$ qui désigne le reste de la division euclidienne de $a$ par
$n$, l'ensemble $\{\modsimple{2^k}{(2n+1)},k\in\N\}$ est fini,
il existe au moins deux entiers distincts $k$ et $l$ tels que $k<l$ et
$\modtext{2^k}{2^l}{(2n+1)}$. Les nombres $2$ et $(2n+1)$ étant premiers entre
eux, $2$ est inversible modulo $(2n+1)$ et on peut donc <<~diviser~>> par $2^k$
pour obtenir $\modtext{2^{l-k}}{1}{(2n+1)}$ avec $l-k\ge 1$.
\fin

\rem
Un théorème d'Euler stipule que, si $\f(2n+1)$ est le nombre d'entiers premiers
avec $2n+1$ compris entre 1 et $2n+1$ (fonction indicatrice d'Euler),
on a $\modtext{2^{\f(2n+1)}}{1}{(2n+1)}$.
Ainsi $\f(2n+1)$ est une période de $f$. Lorsque $2n+1$ est premier,
$\f(2n+1)=2n$ et ceci n'est rien d'autre que le petit théorème de Fermat :
$\modtext{2^{2n}}{1}{(2n+1)}$.

La discussion précédente montre que l'évolution de la $i$\ieme\ carte
($i\in\{1,2,3,\dots,2n\}$) est décrite par l'orbite de $i$ sous l'action de $f$
(ou de manière équivalente de $f^{-1}$) :
\beqa
\cO(i)&=&\{f^k(i),k\in\Z\}
\\
&=&\{i,f(i),f^2(i),\dots,f^{r-1}(i)\}
\\
&=&\{\modsimple{2^ki}{(2n+1)},0\le k\le r-1\}.
\eeqa
L'orbite de $1$ est en particulier
$$
\cO(1)=\{\modsimple{2^k}{(2n+1)},0\le k\le r-1\}.
$$
Par définition de $r$, le cardinal de $\cO(1)$ est exactement $r$ :
$\card\cO(1)=r$. Pour les autres $i$, on a $\card\cO(i)\le r$.
On voit par ailleurs, puisque $\cO(1)\subset\{1,2,3,\dots,2n\}$,
que $r\le 2n$. L'égalité $r=2n$ signifie que l'on a
une seule orbite : pour tout $i\in\{1,2,3,\dots,2n\}$,
$$
\cO(i)=\{1,2,3,\dots,2n\}.
$$
Dans ce cas, une carte donnée visite toutes les places du jeu
avant de revenir à sa position initiale.
\bpr{
\mbox{}\\
\indent $\bullet$
Pour tout $i\in\{1,2,\dots,n\}$, le cardinal de l'orbite de $i$ est un diviseur
de celui de l'orbite de $1$ : $\card\cO(i)$ divise $\card\cO(1)$.
\\
\indent $\bullet$
Dans le cas particulier où $i$ est un nombre premier avec $2n+1$, ces orbites ont même cardinal :
\mbox{$\card\cO(i)=\card\cO(1)$.}
\\
\indent $\bullet$
Si $2n+1$ est premier, toutes les orbites ont même cardinal $r$ et il y a $2n/r$ orbites distinctes.
}
\dem
\\
\indent $\bullet$ Introduisons le pgcd des entiers $i$ et $2n+1$ : $d_i=\mathrm{pgcd}(i,2n+1)$
et posons $r_i=\card\cO(i)$ (on a en particulier $r_1=r$).
Le cardinal $r_i$ est le plus petit entier strictement positif tel que $\modtext{2^ki}{i}{(2n+1)}$.
Cette dernière égalité est équivalente à l'assertion <<~$(2n+1)$ divise $i(2^k-1)$~>>, ou encore
<<~$(2n+1)/d_i$ divise $(i/d_i)(2^k-1)$~>>. Par définition de $d_i$, les entiers $i/d_i$ et $m_i=(2n+1)/d_i$
sont premiers entre eux, ce qui montre que $m_i$ divise $2^k-1$.
On voit ainsi que $r_i=\min\{k\in\N^*: \modtext{2^k}{1}{m_i}\}$.
\\
Prouvons que $r_i$ divise $r$. Comme $m_i$ divise \mbox{$2n+1$}, on a l'implication
$\modtext{x}{y}{(2n+1)}\Longrightarrow\modtext{x}{y}{m_i}$ ; on
peut alors définir de manière cohérente,
entre les groupes multiplicatifs
$$
\cO(1)=\{\modsimple{2^k}{(2n+1)},0\le k\le r-1\}
$$
et
$$
G_i=\{\modsimple{2^k}{m_i},0\le k \le r_i-1\},
$$
un morphisme
$$
\phi_i: \begin{array}[t]{ccc} \cO(1) &\longrightarrow & G_i
\\
\modsimple{2^k}{(2n+1)}& \longmapsto & \modsimple{2^k}{m_i}
\end{array}
$$
qui est clairement surjectif. Les groupes $\cO(1)/\ker\phi_i$ et $G_i$ sont donc
isomorphes ce qui montre bien que $\card(G_i)$ (qui vaut $r_i$)
divise $\card \cO(1)$.
\\
\indent $\bullet$
Lorsque $i$ est premier avec $2n+1$, $i$ est inversible modulo $(2n+1)$
et l'application de multiplication par~$i$
$$
\begin{array}{ccc}
\cO(1) & \longrightarrow & \cO(i)
\\
\modsimple{2^k}{(2n+1)} & \longmapsto & \modsimple{2^ki}{(2n+1)}
\end{array}
$$
est une bijection entre les orbites $\cO(1)$ et $\cO(i)$.
Les orbites de $1$ et de $i$ ont donc même cardinal (qui vaut~$r$).
\\
\indent $\bullet$
Enfin, si $2n+1$ est un nombre premier, tout $i\in\{1,2,\dots,2n\}$ est premier
avec $2n+1$ et donc $r_i=r$. La réunion des différentes orbites coïncidant avec
$\{1,2,\dots,2n\}$, il y a alors $2n/r$ orbites distinctes.
\fin

\exam
\textsl{Cas d'un jeu de $52$ cartes (quatre séries de cartes numérotées de
$2$ à $10$, valet, dame, roi, as, de couleurs carreau, pique, c{\oe}ur, trèfle).}
Ce cas correspond à $n=26$. L'étude de l'in-shuffle de ce jeu s'effectue modulo
$2n+1=53$ qui est premier. En renumérotant toutes les cartes de $1$ à $52$,
il y a une seule orbite : $\cO(1)=\{1,2,3,\dots,52\}$ et l'in-shuffle est de
période $52$.

\exam
\textsl{Cas d'un jeu de $54$ cartes (jeu précédent avec deux jokers différents).}
Ce cas correspond à $n=27$. L'étude de l'in-shuffle de ce mélange se mène modulo $2n+1=55$.
On trouve, en renumérotant toutes les cartes de $1$ à $54$, quatre orbites distinctes :
\beqa
\cO(1)&=&\{1,2,4,7,8,9,13,14,16,17,
\\
&&
18,26,28,31,32,34,36,43,49,52\},
\\
\cO(3)&=&\{3,6,12,19,21,23,24,27,29,37,
\\
&&
38,39,41,42,46,47,48,51,53,54\},
\\
\cO(5)&=&\{5,10,15,20,25,30,35,40,45,50\},
\\
\cO(11)&=&\{11,22,33,44\}.
\eeqa
On vérifie que les cardinaux de $\cO(3),\cO(5),\cO(11)$,
valant respectivement $20,10,4$, divisent le cardinal de $\cO(1)$ qui vaut $20$.
Ainsi, l'in-shuffle est de période $20$.

\bco{\label{cor-in-shuffle}
\mbox{}\\
\indent $\bullet$
La période de l'in-shuffle d'un jeu de $2^p$ cartes ($p\ge 1$) est $2p$.
\\
\indent $\bullet$
La période de l'in-shuffle d'un jeu de $2^p-2$ cartes ($p\ge 2$) est $p$.
}
\dem
\'Ecrivons l'orbite de $1$.
\\
\indent $\bullet$
Si $n=2^{p-1}$, alors $2n+1=2^p+1$. On voit que
$\modtext{2^p}{-1}{(2n+1)}$, puis, en élevant au carré, que
$\modtext{2^{2p}}{1}{(2n+1)}$. D'après le théorème~\ref{th-in-shuffle},
on tire que $2p$ est une période de l'in-shuffle. On a ensuite, pour tout
$k\in\{p,\lb p+1,\dots,2p\}$, $\modtext{2^k}{2^p-2^{k-p}+1}{(2n+1)}$ avec
l'encadrement $1\le 2^p-2^{k-p}+1\le 2n$. L'orbite de $1$ est dans ce cas
\beqa
\cO(1)
&=&
\{1,2,2^2,\dots,2^{p-1},2^p,\\
&&
\modsimple{2^{p+1}}{(2n+1)},\modsimple{2^{p+2}}{(2n+1)},\\
&&
\dots,\modsimple{2^{2p-1}}{(2n+1)}\}
\\
&=&
\{1,2,2^2,\dots,2^{p-1},2^p,
\\
&&
2^p-1,2^p-3,\dots,2^p-2^{p-1}+1\}
\\
&=&
\{2^k,0\le k\le p-1\}\\
&&
\cup\, \{2^p-2^k+1,0\le k\le p-1\}.
\eeqa
On constate que $\card\cO(1)=2p$ et donc que $2p$ est la période de l'in-shuffle.
\\
\indent $\bullet$
Si $n=2^{p-1}-1$, alors $2n+1=2^p-1$ et donc $\modtext{2^p}{1}{(2n+1)}$.
Cette fois, l'orbite de $1$ est donnée par
$$
\cO(1)=\{1,2,2^2,\dots,2^{p-1}\}
$$
et l'on obtient $\card\cO(1)=p$. Le nombre $p$ est la période de l'in-shuffle.
\fin

\rem
D'après la remarque suivant le théorème~\ref{th-in-shuffle},
on a prouvé au passage que $p$ divise $\f(2^p-1)$ et $2p$ divise $\f(2^p+1)$.

\exam
\textsl{Cas d'un jeu de $32$ cartes (quatre séries de cartes numérotées de
$7$ à $10$, valet, dame, roi, as, de couleurs carreau, pique, c{\oe}ur, trèfle).}
Ce cas correspond à $n=16$ et $p=5$. L'étude de l'in-shuffle de ce jeu se réalise
modulo $2n+1=33$. En renumérotant les cartes de $1$ à $32$, l'\'evolution
progressive de la carte \no~$1$ est décrite selon le cycle
$$
\hspace{-.5em}\begin{array}{l}
f^{-1}(1)=2,\,f^{-2}(1)=4,\,f^{-3}(1)=8,\,f^{-4}(1)=16,
\\[1ex]
f^{-5}(1)=32,\,f^{-6}(1)=31,\,f^{-7}(1)=29,
\\[1ex]
f^{-8}(1)=25,\,f^{-9}(1)=17,\,f^{-10}(1)=1.
\end{array}
$$
Plus généralement, on trouve quatre orbites distinctes :
\beqa
\cO(1)&=&\{1,2,4,8,16,17,25,29,31,32\},
\\
\cO(3)&=&\{3,6,9,12,15,18,21,24,27,30\},
\\
\cO(5)&=&\{5,7,10,13,14,19,20,23,26,28\},
\\
\cO(11)&=&\{11,22\}.
\eeqa
On vérifie que les cardinaux de $\cO(3),\cO(5),\cO(11)$, valant respectivement
$10,10,2$, divisent le cardinal de $\cO(1)$ qui vaut $10$. Ainsi, l'in-shuffle
de ce jeu est de période $10$ qui coïncide précisément avec $2p$.

\subsection{\textsf{Cas de l'out-shuffle}}

\bth{\label{th-out-shuffle}
La période de l'out-shuffle d'un jeu de $2n$ cartes est le plus petit entier
$s\ge 1$ vérifiant $\modtext{2^s}{1}{(2n-1)}$. En d'autres termes, $s$ est
l'ordre de $2$ modulo $(2n-1)$.
}
\dem
Rappelons la relation remarquable pour la permutation $\tg$ associée à l'out-shuffle :
$$
\modulo{\tg^{-1}(j)}{2j}{(2n-1)}
$$
qui donne pour tout $k\in\N$,
$$
\modulo{\tg^{-k}(j)}{2^kj}{(2n-1)}.
$$
La période de $\tg$ et donc de $g$ est alors le plus petit entier $s\ge 1$ vérifiant
$\modtext{2^s}{1}{(2n-1)}$.
\fin

\rem
Le théorème~\ref{th-out-shuffle} découle également du
théorème~\ref{th-in-shuffle} : en effet, il a été observé
qu'un out-shuffle de $2n$ cartes est identique à l'in-shuffle du jeu de $2(n-1)$ cartes
obtenu en retirant les première et dernière cartes.

La discussion précédente montre que l'évolution de la $i$\ieme\ carte ($i\in\{0,1,2,\dots,2n-1\}$)
est décrite par l'orbite de $i$ sous l'action de $\tg$ (ou de manière équivalente de $\tg^{-1}$) :
\beqa
\tcO(i)&=&\{\tg^k(i),k\in\Z\}
\\
&=&\{i,\tg(i),\tg^2(i),\dots,\tg^{s-1}(i)\}
\\
&=&\{\modsimple{2^ki}{(2n-1)},0\le k\le s-1\}.
\eeqa
L'orbite de $1$ est en particulier
$$
\tcO(1)=\{\modsimple{2^k}{(2n-1)},0\le k\le s-1\}.
$$

\exam
\textsl{Cas d'un jeu de $52$ cartes.} L'étude de l'out-shuffle de ce jeu
s'effectue modulo \mbox{$2n-1=51$.}
On trouve, en numérotant les cartes de $0$ à $51$, les neufs orbites suivantes :
\beqa
\tcO(0)&=&\{0\},
\\
\tcO(1)&=&\{1,2,4,8,13,16,26,32\},
\\
\tcO(3)&=&\{3,6,12,24,27,39,45,48\},
\\
\tcO(5)&=&\{5,7,10,14,20,28,29,40\},
\\
\tcO(9)&=&\{9,15,18,21,30,33,36,42\},
\\
\tcO(11)&=&\{11,22,23,31,37,41,44,46\},
\\
\tcO(17)&=&\{17,34\},
\\
\tcO(19)&=&\{19,25,35,38,43,47,49,50\},
\\
\tcO(51)&=&\{51\}.
\eeqa
Les cardinaux de ces orbites sont $1$, $2$ ou $8$. La période de l'out-shuffle est $s=8$.

\exam
\textsl{Cas d'un jeu de $54$ cartes.} L'étude de l'out-shuffle de ce jeu
se mène modulo \mbox{$2n-1=53$} qui est premier. On trouve à présent,
en numérotant les cartes de $0$ à $53$, trois orbites distinctes :
$$
\tcO(0)\!=\!\{0\},\,\tcO(1)\!=\!\{1,2,3,\dots,52\},\,\tcO(53)\!=\!\{53\}.
$$
La période dans ce cas est $s=52$.

La remarque suivant le théorème~\ref{th-out-shuffle} fournit l'analogue suivant du corollaire~\ref{cor-in-shuffle}.
\bco{
\mbox{}\\
\indent $\bullet$
La période de l'out-shuffle d'un jeu de $2^p$ cartes ($p\ge 1$) est $p$.
\\
\indent $\bullet$
La période de l'out-shuffle d'un jeu de $2^p+2$ cartes ($p\ge 1$) est $2p$.
}

\exam
\textsl{Cas d'un jeu de $32$ cartes.} Ce cas correspond à $n=16$ et $p=5$.
L'étude de l'out-shuffle se réalise modulo $2n-1=31$.
En numérotant les cartes de $0$ à $31$, l'\'evolution de la carte \no $1$ est décrite selon
le cycle
$\tg^{-1}(1)=2,\,\tg^{-2}(1)=4,\lb\tg^{-3}(1)=8,\,\tg^{-4}(1)=16,\,\tg^{-5}(1)=1.$
Plus généralement, on trouve huit orbites distinctes :
\beqa
\tcO(0)&=&\{0\},
\\
\tcO(1)&=&\{1,2,4,8,16\},
\\
\tcO(3)&=&\{3,6,12,17,24\},
\\
\tcO(5)&=&\{5,9,10,18,20\},
\\
\tcO(7)&=&\{7,14,19,25,28\},
\\
\tcO(11)&=&\{11,13,21,22,26\},
\\
\tcO(15)&=&\{15,23,27,29,30\},
\\
\tcO(31)&=&\{31\}.
\eeqa
Les cardinaux de ces orbites valent $1$ ou $5$ et
l'out-shuffle de ce jeu est de période $5$ qui coïncide précisément avec $p$.

\subsection{\textsf{Cas d'un mélange de Monge}}

\bth{\label{th-monge-shuffle1}
La période commune de $h_1$ et de $h_3$ est le plus petit entier $u\ge 1$ vérifiant
l'une des deux congruences
$\modtext{2^u}{1}{(4n+1)}$ ou $\modtext{2^u}{-1}{(4n+1)}$.
}
\noindent\textsc{Explication du théorème.}
Ou bien il existe une puissance de $2$ congrue à $-1$ modulo \mbox{$(4n+1)$}
et l'on choisit pour $u$ la plus petite. Dans ce cas, $2u$ est l'ordre de
$2$ modulo $(4n+1)$ : $\modtext{2^{2u}}{1}{(4n+1)}$.
Ou bien il n'existe aucune puissance de $2$ congrue à $-1$ modulo $(4n+1)$
et l'on choisit alors pour $u$ l'ordre de $2$ modulo $(4n+1)$.

\dem
Rappelons tout d'abord que $h_1$ et $h_3$ sont conjuguées :
$h_3=s\circ h_1\circ s^{-1}$. Cette propriété s'étend à toutes les itérées :
\mbox{$h^k_3=s\circ h^k_1\circ s^{-1}$.}
Ainsi l'équation $h_1^k=id$ est équivalente à l'équation $h_3^k=id$, prouvant
que $h_1$ et $h_3$ ont même période.
Rappelons la relation
$$
\modulo{h_3^{-1}(j)}{\pm 2j}{(4n+1)}.
\vspace{-.1\baselineskip}
$$
On a
\vspace{-.1\baselineskip}
$$
\modulo{h_3^{-k}(j)}{\pm 2^kj}{(4n+1)},
$$
le signe $\pm$ dépendant de $j$ et $k$ est déterminé sans ambiguïté
par la contrainte $1\le h_3^{-k}(j)\le 2n$.
Plus précisément, si $2^kj$ admet un représentant modulo $(4n+1)$ compris entre
$1$ et $2n$, alors $h_3^{-k}(j)$ coïncide avec ce représentant
et on choisit le signe $+$ dans la congruence : $\modtext{h_3^{-k}(j)}{+2^kj}{(4n+1)}$.
Sinon, $2^kj$ admet un représentant entre $2n+1$ et $4n$ ($0$ ne peut pas être
un représentant car $2^k$ et $(4n+1)$ sont premiers entre eux et $1\le j\le 2n$).
Dans ce cas, $-2^kj$ admet un représentant entre $1$ et $2n$ et l'on choisit le
signe $-$ : $\modtext{h_3^{-k}(j)}{-2^kj}{(4n+1)}$.
\\
\indent
Pour avoir en particulier $h_3^{-k}(1)=1$,
on doit choisir $k$ tel que $\modtext{2^k}{\pm 1}{(4n+1)}$.
Considérons donc le plus petit entier $u\ge 1$ tel que
\mbox{$\modtext{2^u}{1}{(4n+1)}$ ou $\modtext{2^u}{\!-1}{\!(4n+1)}$.}
On a ensuite pour tout $j\in\{1,2,\dots,2n\}$,
$\modtext{h_3^{-u}(j)}{\pm j}{(4n+1)}$,
le signe $\pm$ dépendant de $j$ et $u$. En fait, le représentant de $\modsimple{-j}{(4n+1)}$
compris entre $1$ et $4n+1$ est $4n+1-j$. Mais ce représentant est supérieur
à $2n$ et la condition $1\le h_3^{-u}(j)\le 2n$ n'est remplie que dans le cas où
$h_3^{-u}(j)=j$.
Ainsi $h_3^{-u}=id$ et $u$ est la période de $h_3$.
\fin

\exam
\textsl{Cas d'un jeu de $52$ cartes.}
L'étude de ce mélange de Monge s'effectue modulo \mbox{$4n+1=105$.}
On trouve les orbites suivantes :
\beqa
\cO(1)&=&\{1,2,4,8,13,16,23,26,32,41,46,52\},
\\
\cO(3)&=&\{3,6,9,12,18,24,27,33,36,39,48,51\},
\\
\cO(5)&=&\{5,10,20,25,40,50\},
\\
\cO(7)&=&\{7,14,28,49\},
\\
\cO(11)&=&\{11,17,19,22,29,31,
\\
&&34,37,38,43,44,47\},
\\
\cO(15)&=&\{15,30,45\},
\\
\cO(21)&=&\{21,42\},
\\
\cO(35)&=&\{35\}.
\eeqa
Ce mélange a pour période $12$.

\exam
\textsl{Cas d'un jeu de $54$ cartes.}
L'étude de ce mélange de Monge s'effectue modulo \mbox{$4n+1=109$.}
On trouve les trois orbites suivantes :
\beqa
\cO(1)&=&\{1,2,4,8,16,17,19,23,27,
\\
&&
32,33,34,38,41,43,45,46,54\},
\\
\cO(3)&=&\{3,5,6,7,10,12,13,14,20,
\\
&&
24,26,28,29,40,48,51,52,53\},
\\
\cO(9)&=&\{9,11,15,18,21,22,25,30,31,
\\
&&
35,36,37,39,42,44,47,49,50\},
\eeqa
Ce mélange a pour période $18$.

De manière similaire, la congruence précédemment signalée
$$
\modulo{\th_4^{-1}(j)}{\pm 2j}{(4n-1)}
$$
fournit la période de $h_4$.
\bth{\label{th-monge-shuffle2}
La période commune de $h_2$ et de $h_4$ est le plus petit entier $v\ge 1$ vérifiant
l'une des deux congruences $\modtext{2^v}{1}{(4n-1)}$ ou $\modtext{2^v}{-1}{(4n-1)}$.
}

\exam
\textsl{Cas d'un jeu de $52$ cartes.}
L'étude de ce mélange de Monge s'effectue modulo \mbox{$4n-1=103$.}
On trouve les deux orbites \mbox{$\tcO(0)=\{0\}$} et $\tcO(1)=\{1,2,3,\dots,51\}$.
La~période de ce mélange est $51$.

\exam
\textsl{Cas d'un jeu de $54$ cartes.}
L'étude de ce mélange de Monge s'effectue modulo \mbox{$4n-1=107$.}
On trouve les deux orbites \mbox{$\tcO(0)=\{0\}$} et $\tcO(1)=\{1,2,3,\dots,53\}$.
La~période de ce mélange est $53$.

\bco{
La période d'un mélange de Monge de $2^p$ cartes ($p\ge 1$) est $p+1$.
}
\dem
Pour $n=2^{p-1}$, on a $4n+1\linebreak =2^{p+1}+1$. Donc $\modtext{2^{p+1}}{-1}{(4n+1)}$
et $\modtext{2^{p+1}}{1}{(4n-1)}$.
D'après les théorèmes~\ref{th-monge-shuffle1} et~\ref{th-monge-shuffle2}, $p+1$ est une
période de $h_1,h_2,h_3$ et $h_4$.
Dans les deux cas, l'orbite de $1$ est l'ensemble $\{1,2,2^2,\dots,2^p\}$ de cardinal $p+1$,
ce qui prouve le corollaire.
\fin

\exam
\textsl{Cas d'un jeu de $32$ cartes.}
L'étude du mélange de Monge associé à $h_3$ s'effectue modulo $4n+1=65$.
On trouve les orbites suivantes :
\beqa
\cO(1)&=&\{1,2,4,8,16,32\},
\\
\cO(3)&=&\{3,6,12,17,24,31\},
\\
\cO(5)&=&\{5,10,15,20,25,30\},
\\
\cO(7)&=&\{7,9,14,18,28,29\},
\\
\cO(11)&=&\{11,19,21,22,23,27\},
\\
\cO(13)&=&\{13,26\}.
\eeqa
L'étude du mélange de Monge associé à $\th_4$ s'effectue modulo $4n-1=63$.
On trouve les orbites suivantes :
\beqa
\tcO(0)&=&\{0\},
\\
\tcO(1)&=&\{1,2,4,8,16,31\},
\\
\tcO(3)&=&\{3,6,12,15,24,30\},
\\
\tcO(5)&=&\{5,10,17,20,23,29\},
\\
\tcO(7)&=&\{7,14,28\},
\\
\tcO(9)&=&\{9,18,27\},
\\
\tcO(11)&=&\{11,13,19,22,25,26\},
\\
\tcO(21)&=&\{21\}.
\eeqa
Dans les deux cas, la période est $6$ et coïncide bien avec $p+1$
(ici $n=2^{p-1}$ avec $p=5$).

\subsection{\textsf{Quelques valeurs numériques}}\label{section-num}

Nous donnons ci-dessous les périodes des in-shuffles et des mélanges de Monge
associés à $h_1$ et $h_2$ (baptisés dans le tableau de <<~in-Monge~>> et
<<~out-Monge~>>) pour des jeux de $2n$ cartes avec $2n\le 64$.

\vspace{\baselineskip}
\newfont{\mafonte}{cmr10}
\vbox{\mafonte
\noindent\begin{tabular}{l|cccccccc}
\hline
$2n$ &2&4&6&8&10&12&14&16
\\
\hline\hline
in-shuffle &2&4&3&6&10&12&4&8
\\
\hline
in-Monge &2&3&6&4&6&10&14&5
\\
\hline
out-Monge &1&3&5&4&9&11&9&5
\\
\hline
\end{tabular}
}

\vspace{\baselineskip}
\vbox{\mafonte
\noindent\begin{tabular}{l|cccccccc}
\hline
$2n$ &18&20&22&24&26&28&30&32
\\
\hline\hline
in-shuffle &18&6&11&20&18&28&5&10
\\
\hline
in-Monge &18&10&12&21&26&9&30&6
\\
\hline
out-Monge &12&12&7&23&8&20&29&6
\\
\hline
\end{tabular}
}

\vspace{\baselineskip}
\vbox{\mafonte
\noindent\begin{tabular}{l|cccccccc}
\hline
$2n$ &34&36&38&40&42&44&46&48
\\
\hline\hline
in-shuffle &12&36&12&20&14&12&23&21
\\
\hline
in-Monge &22&9&30&27&8&11&10&24
\\
\hline
out-Monge &33&35&20&39&41&28&12&36
\\
\hline
\end{tabular}
}

\vspace{\baselineskip}
\vbox{\mafonte
\noindent\begin{tabular}{l|cccccccc}
\hline
$2n$ &50&52&54&56&58&60&62&64
\\
\hline\hline
in-shuffle &8&52&20&18&58&60&6&12
\\
\hline
in-Monge &50&12&18&14&12&55&50&7
\\
\hline
out-Monge &15&51&53&36&44&24&20&7
\\
\hline
\end{tabular}
}

\section{\textsf{In-shuffle : cas particuliers}}\label{section-caspart1}

Dans cette partie, nous calculons explicitement les itérations successives de
la permutation associée à l'in-shuffle dans les deux cas particuliers
$n=2^{p-1}$ et $n=2^{p-1}-1$. L'astuce de calcul consiste à travailler
avec les écritures binaires des numéros de cartes.

Nous ne calculerons pas les itérations successives de l'out-shuffle dans les
cas $n=2^{p-1}$ et $n=2^{p-1}+1$, l'out-shuffle étant équivalent à l'in-shuffle
associé aux cas respectifs $n=2^{p-1}-1$ et $n=2^{p-1}$.

\subsection{\textsf{Cas d'un jeu de $2^p$ cartes}}

Nous nous plaçons dans le cas d'un jeu de $2^p$ cartes, c'est-à-dire $n=2^{p-1}$.
Nous travaillons ici avec $\tf$.
Introduisons l'écriture binaire d'un \mbox{$i\in\{0,1,\dots,2n-1\}$ :}
$$
i=\overline{i_{p-1}\dots i_0}=\sum_{k=0}^{p-1} i_k2^k
$$
où les $i_0,i_1,\dots,i_{p-1}$ sont des bits $0$ ou $1$. On a en particulier
$n=\overline{1\underset{p-1}{\underbrace{0-0}}}$
et
$2n-1=\overline{\underset{p}{\underbrace{1-1}}}$.
Dans ces conditions, l'image de $i$ par la permutation $\tf$ s'écrit
$$
\tf(i)=\left\{\begin{array}{ll}
\overline{1i_{p-1}\dots i_1} & \mbox{si $i_0=0$,}
\\[1ex]
\overline{i_{p-1}\dots i_1} & \mbox{si $i_0=1$,}
\end{array}\right.
$$
soit encore
$$
\tf(i)=\overline{(1-i_0)i_{p-1}\dots i_1}.
$$

\subsubsection{\textsf{Calcul des itérées $\tf^k(i)$ pour chaque $i$}}

Avec cette représentation de $\tf$, on voit progressivement que
\beqa
\tf(i)&=&\overline{(1-i_0)i_{p-1}\dots i_1},
\\
\tf^2(i)&=&\overline{(1-i_1)(1-i_0)i_{p-1}\dots i_2}
\eeqa
et plus généralement, pour $0\le k\le p$,
$$
\tf^k(i)=\overline{(1-i_{k-1})\dots(1-i_0)i_{p-1}\dots i_k}.
$$
Pour $k=p$, on obtient
$$
\tf^p(i)=\overline{(1-i_{p-1})\dots (1-i_0)}
$$
que l'on peut réécrire
$$
\tf^p(i)=\overline{\underset{p}{\underbrace{1-1}}}-\overline{i_{p-1}\dots i_0}
=2n-1-i.
$$
Ainsi $\tf^p$ est la symétrie des entiers $0,1,\dots,2n-1$.
Cela signifie que dans le cas d'une pile de $2n$ jetons dont les $n$ jetons du
bas sont verts et les $n$ du haut sont rouges, on obtient
au bout de $p$ mélanges parfaits une pile inversée : les $n$ jetons
du bas sont rouges et les $n$ du haut sont verts.

En poursuivant, on trouve pour $p\le k\le 2p$,
$$
\tf^k(i)=\overline{i_{k-p-1} \dots i_0(1-i_{p-1})\dots (1-i_{k-p})}
$$
pour aboutir finalement à
$$
\tf^{2p}(i)=\overline{i_{p-1} i_{p-2}\dots i_0}=i.
$$
On retrouve bien le fait que, dans le cas où \mbox{$n=2^{p-1}$,} le nombre $2p$ est une
période de $\tf$ (et aussi de $f$).

\subsubsection{\textsf{Calcul des itérées $\tf^k(0)$}}

L'évolution en binaire de la première carte (numéro $0$) est intéressante en soi.
Les calculs précédents donnent
$$
\tf(0)=\overline{1\underset{p-1}{\underbrace{0-0}}}=n.
$$
Plus généralement, pour $0\le k\le p$,
\beqa
\tf^k(0)&=&\overline{\underset{k}{\underbrace{1-1}}\underset{p-k}{\underbrace{0-0}}}
\\
&=&2^{p-1}+2^{p-2}+\dots+2^{p-k}=2^p-2^{p-k}
\eeqa
qui conduit à
$$
\tf^p(0)=\overline{\underset{p}{\underbrace{1-1}}}=2^p-1=2n-1.
$$
Puis, pour $p\le k\le 2p$,
$$
\tf^k(0)=\overline{\underset{2p-k}{\underbrace{1-1}}}=2^{2p-k}-1
$$
et finalement $\tf^{2p}(0)=0$.
Ainsi, l'orbite de $0$ sous l'action de $\tf$ est
\beqa
\tcO(0)&=&\{2^p-2^{p-k},1\le k \le p\}
\\
&&
\cup\, \{2^{2p-k}-1,p\le k\le 2p-1\}
\\
&=&\{2^k-1,1\le k \le p\}
\\
&&
\cup\, \{2^p-2^k,1\le k\le p\}.
\eeqa
On voit que $\card\tcO(0)=2p$. Le nombre $2p$ est bien la période de $\tf$.

\subsection{\textsf{Cas d'un jeu de $2^p-2$ cartes}}

Nous nous plaçons dans le cas où $n=2^{p-1}-1$.
Nous travaillons ici avec $f$.
Introduisons de nouveau l'écriture binaire d'un $i\in\{1,2\dots,2n\}$ : $i=\overline{i_{p-1}\dots i_0}$.
On a $n+1=2^{p-1}=\overline{1\underset{p-1}{\underbrace{0-0}}}$.
Dans ces conditions, l'image de $i$ par la permutation $f$ s'écrit
$$
f(i)=\left\{\begin{array}{ll}
\overline{i_{p-1}\dots i_1} & \mbox{si $i_0=0$,}
\\[1ex]
\overline{1i_{p-1}\dots i_1} & \mbox{si $i_0=1$,}
\end{array}\right.
$$
soit encore
$$
f(i)=\overline{i_0i_{p-1}\dots i_1}.
$$
Dans ce cas, la permutation $f$ correspond à une simple permutation circulaire
des chiffres de la décomposition binaire de la variable :
$(i_0,i_1,\dots,i_{p-1})\longmapsto (i_1,i_2\dots,i_{p-1},i_0)$.

\subsubsection{\textsf{Calcul des itérées $f^k(i)$ pour chaque $i$}}

On obtient immédiatement, pour $0\le k\le p$,
$$
f^k(i)=\overline{i_{k-1}i_{k-2}\dots i_0i_{p-1}\dots i_k}.
$$
Finalement pour $k=p$ :
$$
f^p(i)=\overline{i_{p-1}i_{p-2}\dots i_0}=i.
$$
Ce procédé est signalé dans \cite{biane,diaconis-graham-kantor}. On retrouve
le fait que, dans le cas où $n=2^{p-1}-1$, $p$ est une période de $f$.

\subsubsection{\textsf{Calcul des itérées $f^k(1)$}}

Examinons l'évolution binaire de la première carte (\no $1$) : la permutation circulaire
des chiffres donne
$$
f(1)=\overline{1\underset{p-1}{\underbrace{0-0}}}=2^{p-1}=n+1,
$$
puis, pour $0\le k\le p$,
$$
f^k(1)=\overline{1\underset{p-k}{\underbrace{0-0}}}=2^{p-k}.
$$
On arrive finalement à $f^p(1)=1$ et $p$ est bien la période de $f$.

\section{\textsf{Mélange de Monge : cas d'un jeu de $2^p$ cartes}}\label{section-caspart2}

Dans toute cette section, nous nous plaçons dans le cas où $n=2^{p-1}$.

\subsection{\textsf{Version associée à $h_1$}}

Nous travaillons ici avec la permutation $\th_1$ des entiers $0,1,\dots,2n-1$.
On a, pour $i=\overline{i_{p-1}\dots i_0}$,
$$
\th_1(i)=\!\left\{\!\!\begin{array}{ll}
\overline{(1-i_0)i_{p-1}\dots i_1} & \mbox{si $i_0=0$,}
\\[1ex]
\overline{(1-i_0)(1-i_{p-1})\dots (1-i_1)} & \mbox{si $i_0=1$.}
\end{array}\right.
$$
L'expression de $\th_1(i)$ ci-dessus dépendant de la parité de $i$,
nous sommes amenés à décomposer l'écriture binaire de $i$ en blocs
de $0$ et de $1$ consécutifs.
Notons $m$ ($m\ge 1$) le nombre de blocs apparaissant dans $i$ et $l_1,l_2,\dots,l_m$
leurs longueurs respectives en partant de la droite vers la gauche
($l_1,\dots,l_m\ge 1$ et $l_1+\dots+l_m=p$). Pour $m=1$, on a les deux possibilités
$i=\overline{0-0}$ et $i=\overline{1-1}$.
Pour $m\ge 2$, on a les quatre possibilités génériques :
\beqa
i&=&\overline{\underset{l_m}{\underbrace{0-0}}\underset{l_{m-1}}{\underbrace{1-1}}
\dots\underset{l_2}{\underbrace{1-1}}\underset{l_1}{\underbrace{0-0}}},
\\
i&=&\overline{\underset{l_m}{\underbrace{0-0}}\underset{l_{m-1}}{\underbrace{1-1}}
\dots\underset{l_2}{\underbrace{0-0}}\underset{l_1}{\underbrace{1-1}}},
\\
i&=&\overline{\underset{l_m}{\underbrace{1-1}}\underset{l_{m-1}}{\underbrace{0-0}}
\dots\underset{l_2}{\underbrace{1-1}}\underset{l_1}{\underbrace{0-0}}},
\\
i&=&\overline{\underset{l_m}{\underbrace{1-1}}\underset{l_{m-1}}{\underbrace{0-0}}
\dots\underset{l_2}{\underbrace{0-0}}\underset{l_1}{\underbrace{1-1}}}.
\\
\eeqa

\subsubsection{\textsf{Calcul des itérées $\th_1^k(i)$ pour chaque $i$}}

Nous considérons par exemple le cas d'un nombre $i$ se décomposant sous la forme
$$
i=\overline{\underset{l_m}{\underbrace{1-1}}\underset{l_{m-1}}{\underbrace{0-0}}
\dots\underset{l_2}{\underbrace{1-1}}\underset{l_1}{\underbrace{0-0}}}.
$$
On a
\beqa
\th_1(i)&=&\overline{\underset{l_m+1}{\underbrace{1-1}}\underset{l_{m-1}}{\underbrace{0-0}}
\dots\underset{l_2}{\underbrace{1-1}}\underset{l_1-1}{\underbrace{0-0}}}
\\[-4ex]
&\vdots&
\\
\th_1^{l_1}(i)&=&\overline{\underset{l_m+l_1}{\underbrace{1-1}}\underset{l_{m-1}}{\underbrace{0-0}}
\dots\underset{l_2}{\underbrace{1-1}}}
\eeqa
puis
\beqa
\th_1^{l_1+1}(i)&=&\overline{\underset{l_m+l_1+1}{\underbrace{0-0}}\underset{l_{m-1}}{\underbrace{1-1}}
\dots\underset{l_2-1}{\underbrace{0-0}}}
\\
\th_1^{l_1+2}(i)&=&\overline{1\underset{l_m+l_1+1}{\underbrace{0-0}}\underset{l_{m-1}}{\underbrace{1-1}}
\dots\underset{l_2-2}{\underbrace{0-0}}}
\\[-4ex]
&\vdots&
\\
\th_1^{l_1+l_2}(i)&=&\overline{\underset{l_2-1}{\underbrace{1-1}}\underset{l_m+l_1+1}{\underbrace{0-0}}
\underset{l_{m-1}}{\underbrace{1-1}}\dots\underset{l_3}{\underbrace{1-1}}}
\eeqa
puis
\beqa
\th_1^{l_1+l_2+1}(i)&=&\overline{\underset{l_2}{\underbrace{0-0}}\underset{l_m+l_1+1}{\underbrace{1-1}}
\underset{l_{m-1}}{\underbrace{0-0}}\dots\underset{l_3-1}{\underbrace{0-0}}}
\\[-4ex]
&\vdots&
\\
\th_1^{l_1+l_2+l_3}(i)&=&\overline{\underset{l_3-1}{\underbrace{1-1}}
\underset{l_2}{\underbrace{0-0}}\underset{l_m+l_1+1}{\underbrace{1-1}}
\underset{l_{m-1}}{\underbrace{0-0}}\dots\underset{l_4}{\underbrace{1-1}}}.
\eeqa
On continue ainsi de suite pour arriver à
\beqa
\th_1^{l_1+\dots+l_{m-1}}(i)&=&\overline{\underset{l_{m-1}-1}{\underbrace{1-1}}
\underset{l_{m-2}}{\underbrace{0-0}}\dots\underset{l_2}{\underbrace{0-0}}\underset{l_m+l_1+1}{\underbrace{1-1}}}
\\
\th_1^{l_1+\dots+l_{m-1}+1}(i)&=&\overline{\underset{l_{m-1}}{\underbrace{0-0}}
\underset{l_{m-2}}{\underbrace{1-1}}\dots\underset{l_2}{\underbrace{1-1}}\underset{l_m+l_1}{\underbrace{0-0}}}
\\[-4ex]
&\vdots&
\\
\th_1^{l_1+\dots+l_m+1}(i)&=&\overline{\underset{l_m}{\underbrace{1-1}}
\underset{l_{m-1}}{\underbrace{0-0}}\dots\underset{l_2}{\underbrace{1-1}}\underset{l_1}{\underbrace{0-0}}}.
\eeqa
Finalement, le résultat de la dernière étape s'écrit exactement
$$
\th_1^{p+1}(i)=i.
$$
Les calculs menés ci-dessus s'étendent aisément aux autres formes possibles
de décompositions binaires par blocs de $i$ quitte à faire $l_1=0$ ou/et $l_m=0$.
\subsubsection{\textsf{Calcul des itérées $\th_1^k(0)$}}

Les calculs de la section précédente donnent en particulier les itérations
successives de $\th_1$ en $0$ :
$$
\th_1(0)=\overline{1\underset{p-1}{\underbrace{0-0}}}=n,\quad
\th_1^2(0)=\overline{11\underset{p-2}{\underbrace{0-0}}}.
$$
Plus généralement, on a pour $0\le k\le p$,
$$
\th_1^k(0)=\overline{\underset{k}{\underbrace{1-1}}\underset{p-k}{\underbrace{0-0}}}
$$
jusqu'à
$$
\th_1^p(0)=\overline{\underset{p}{\underbrace{1-1}}}=2n-1
$$
et enfin $\th_1^{p+1}(0)=0$.

\subsection{\textsf{Version associée à $h_2$}}

En ce qui concerne la permutation $\th_2$ des entiers $0,1,\dots,2n-1,$ on a
pour $i=\overline{i_{p-1}\dots i_0}$,
$$
\th_2(i)=\left\{\begin{array}{ll}
\overline{i_0(1-i_{p-1})\dots (1-i_1)} & \mbox{si $i_0=0$,}
\\[1ex]
\overline{i_0i_{p-1}\dots i_1} & \mbox{si $i_0=1$,}
\end{array}\right.
$$

\subsubsection{\textsf{Calcul des itérées $\th_2^k(i)$ pour chaque $i$}}

Considérons de nouveau l'exemple où
$$
i=\overline{\underset{l_m}{\underbrace{1-1}}\underset{l_{m-1}}{\underbrace{0-0}}
\dots\underset{l_2}{\underbrace{1-1}}\underset{l_1}{\underbrace{0-0}}}.
$$
On a
\beqa
\th_2(i)&=&\overline{\underset{l_m+1}{\underbrace{0-0}}\underset{l_{m-1}}{\underbrace{1-1}}
\dots\underset{l_2}{\underbrace{0-0}}\underset{l_1-1}{\underbrace{1-1}}}
\\
\th_2^2(i)&=&\overline{1\underset{l_m+1}{\underbrace{0-0}}\underset{l_{m-1}}{\underbrace{1-1}}
\dots\underset{l_2}{\underbrace{0-0}}\underset{l_1-2}{\underbrace{1-1}}}
\\[-4ex]
&\vdots&
\\
\th_2^{l_1}(i)&=&\overline{\underset{l_1-1}{\underbrace{1-1}}\underset{l_m+1}{\underbrace{0-0}}
\underset{l_{m-1}}{\underbrace{1-1}}\dots\underset{l_2}{\underbrace{0-0}}}
\eeqa
puis
\beqa
\th_2^{l_1+1}(i)&=&\overline{\underset{l_1}{\underbrace{0-0}}\underset{l_m+1}{\underbrace{1-1}}
\underset{l_{m-1}}{\underbrace{0-0}}\dots\underset{l_2-1}{\underbrace{1-1}}}
\\[-4ex]
&\vdots&
\\
\th_2^{l_1+l_2}(i)&=&\overline{\underset{l_2-1}{\underbrace{1-1}}\underset{l_1}{\underbrace{0-0}}
\underset{l_m+1}{\underbrace{1-1}}\underset{l_{m-1}}{\underbrace{0-0}}\dots\underset{l_3}{\underbrace{0-0}}}
\eeqa
puis
\beqa
\th_2^{l_1+l_2+1}(i)\!\!&\!\!=\!\!&\!\overline{\underset{l_2}{\underbrace{0-0}}\underset{l_1}{\underbrace{1-1}}
\underset{l_m+1}{\underbrace{0-0}}\underset{l_{m-1}}{\underbrace{1-1}}
\dots\underset{l_3-1}{\underbrace{1-1}}}
\\[-4ex]
\!\!&\!\!\vdots\!\!&\!
\\
\th_2^{l_1+l_2+l_3}(i)\!\!&\!\!=\!\!&\!\overline{\underset{l_3-1}{\underbrace{1-1}}
\underset{l_2}{\underbrace{0-0}}\underset{l_1}{\underbrace{1-1}}
\underset{l_m+1}{\underbrace{0-0}}\underset{l_{m-1}}{\underbrace{1-1}}
\dots\underset{l_4}{\underbrace{0-0}}}.
\eeqa
On continue ainsi de suite pour arriver à
\beqa
\th_2^{l_1+\dots+l_{m-1}}(i)&=&\overline{\underset{l_{m-1}-1}{\underbrace{1-1}}
\underset{l_{m-2}}{\underbrace{0-0}}\dots\underset{l_1}{\underbrace{1-1}}\underset{l_m+1}{\underbrace{0-0}}}
\\
\th_2^{l_1+\dots+l_{m-1}+1}(i)&=&\overline{\underset{l_{m-1}}{\underbrace{0-0}}
\underset{l_{m-2}}{\underbrace{1-1}}\dots\underset{l_1}{\underbrace{0-0}}\underset{l_m}{\underbrace{1-1}}}
\\[-4ex]
&\vdots&
\\
\th_2^{l_1+\dots+l_m}(i)&=&\overline{\underset{l_m-1}{\underbrace{1-1}}
\underset{l_{m-1}}{\underbrace{0-0}}\dots\underset{l_1}{\underbrace{0-0}}1}
\\
\th_2^{l_1+\dots+l_m+1}(i)&=&\overline{\underset{l_m}{\underbrace{1-1}}
\underset{l_{m-1}}{\underbrace{0-0}}\dots\underset{l_2}{\underbrace{1-1}}\underset{l_1}{\underbrace{0-0}}}.
\eeqa
Finalement, la dernière étape donne exactement
$$
\th_2^{p+1}(i)=i.
$$

\subsubsection{\textsf{Calcul des itérées $\th_2^k(0)$}}

Les itérations successives de $\th_2$ en $0$ s'écrivent
$$
\th_2(0)=\overline{0\underset{p-1}{\underbrace{1-1}}}=n-1,\quad
\th_2^2(0)=\overline{10\underset{p-2}{\underbrace{1-1}}}.
$$
Plus généralement, on a pour $0\le k\le p$,
$$
\th_2^k(0)=\overline{\underset{k-1}{\underbrace{1-1}}0\underset{p-k}{\underbrace{1-1}}}
$$
jusqu'à
$$
\th_2^p(0)=\overline{\underset{p-1}{\underbrace{1-1}}0}=2n-2
$$
et enfin $\th_2^{p+1}(0)=0$.

\section{\textsf{Cas d'un jeu contenant un nombre impair de cartes}}\label{section-impair}

On considère ici rapidement le cas d'un jeu de $2n+1$ cartes, cas étudié
dans \cite{morris1}. Dans cette situation, le coupage de ce jeu en deux
parties peut se faire de deux manières :
soit à la $n\ieme$ carte, soit à la $(n+1)\ieme$. En d'autres termes, après
coupage, le premier paquet contient $n$ cartes et le deuxième en contient $n+1$,
ou bien le premier paquet contient $n+1$ cartes et le deuxième en contient $n$.
On intercale alors le paquet de $n$ cartes dans celui de $n+1$.

Examinons les deux situations possibles.

\indent $\bullet$ \textsl{Premier coupage :}
le jeu de cartes numérotées dans l'ordre $1,2,$ $3,\dots,2n,2n+1$ est coupé en
deux paquets de cartes numérotées $1,2,\dots,n$ pour le premier et
$n+1,n+2,\dots,2n,2n+1$ pour le deuxième.
On constitue un nouveau jeu de $2n+1$ cartes en intercalant le premier paquet
dans le deuxième, donnant ainsi la suite de cartes
\nos\ $n+1,1,n+2,2,\dots,n-1,2n,n,2n+1$ (voir Fig.~\ref{fig-impair}).
On observe que la carte \no $(2n+1)$ reste immobile et qu'en la retirant du jeu,
cette manipulation est identique à l'in-shuffle du jeu des $2n$ cartes restantes.
\\
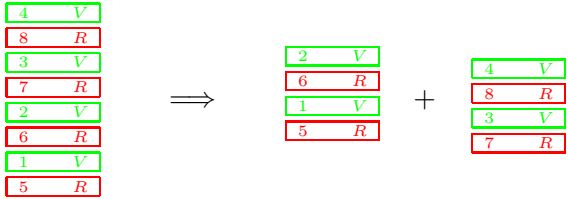
\indent $\bullet$ \textsl{Deuxième coupage :}
le jeu de cartes renumérotées dans l'ordre $0,1,2,3,\dots,2n$ est coupé à
présent en deux paquets de cartes numérotées $0,1,2,\dots,n$ pour le premier et
\mbox{$n+1,n+2,\dots,2n$} pour le deuxième.
On forme un nouveau jeu de $2n+1$ cartes en intercalant le deuxième paquet dans
le premier, fournissant ainsi la suite de cartes
\nos\ $0,n+1,1,n+2,2,\dots,2n-1,\linebreak n-1,2n,n$ (voir Fig.~\ref{fig-impair}).
Dans ce cas, la carte \no $0$ reste immobile ; en la retirant du jeu,
cette manipulation est identique à l'in-shuffle du jeu des $2n$ cartes restantes.

\dessinshuffleimpair

En conclusion, on a le résultat suivant.
\bth{
La période du mélange parfait de $2n+1$ cartes est l'ordre de $2$ modulo $(2n+1)$.
}

\section{\textsf{\mbox{Déplacement d'une carte vers} une position donnée dans
le cas d'un jeu de $2^p$ cartes}}
\label{section-deplacement}

Dans cette section, nous considérons le problème d'Elmsley consistant à
déterminer une succession d'in- et d'out-shuffles déplaçant une carte donnée
à une position donnée. Nous nous plaçons dans le cas simple d'un jeu de $2^p$
cartes ($p\ge 1$) et renvoyons le lecteur à~\cite{diaconis-graham} où une procédure
algorithmique est proposée dans le cas général.
Les cartes sont numérotées de bas en haut $0,1,2,\dots,2^p-1$.

\subsection{\textsf{Procédure générale}}

Rappelons que les déplacements
de la carte \no $j$ où $j=\overline{j_{p-1}\dots j_0}$,
par un in- et un out-shuffles sont représentés, en posant pour simplifier
les notations $\tf^{-1}=\ff$ et $\tg^{-1}=\gg$, par
\beqa
\ff(j)&=&\overline{j_{p-2}\dots j_0(1-j_{p-1})},
\\
\gg(j)&=&\overline{j_{p-2}\dots j_0j_{p-1}}.
\eeqa
Si l'on effectue consécutivement $k_1$ in-, $k_2$ out-, $k_3$ in-, $k_4$ out-,
$\dots$, $k_{m-1}$ in- et $k_m$ out-shuffles où $k_1,k_2,\dots,k_m$ sont des
nombres positifs (éventuellement $k_1,k_m$ pouvant être nuls) de somme $p$,
la carte \no $j$ se retrouve à la position de numéro
$$
i=(\gg^{k_m}\circ\ff^{k_{m-1}}\circ\dots
\circ\gg^{k_4}\circ\ff^{k_3}\circ\gg^{k_2}\circ\ff^{k_1})(j).
$$
Déterminons explicitement l'écriture binaire du numéro $i$. On a
$$
\ff^{k_1}(j)=\overline{j_{p-k_1-1}\dots j_0(1-j_{p-1})\dots(1-j_{p-k_1})}
$$
puis
\beqa
\lqn{(\gg^{k_2}\circ\ff^{k_1})(j)}
\\[-2ex]
&=&\overline{j_{p-k_1-k_2-1}\dots j_0}
\\
&&\overline{(1-j_{p-1})\dots(1-j_{p-k_1})j_{p-k_1-1}\dots j_{p-k_1-k_2}}
\eeqa
puis
\beqa
\lqn{(\ff^{k_3}\circ\gg^{k_2}\circ\ff^{k_1})(j)}
\\[-2ex]
&=&\overline{j_{p-k_1-k_2-k_3-1}\dots j_0}
\\
&&\overline{(1-j_{p-1})\dots(1-j_{p-k_1})}
\\
&&\overline{j_{p-k_1-1}\dots j_{p-k_1-k_2}}
\\
&&\overline{(1-j_{p-k_1-k_2-1})\dots(1-j_{p-k_1-k_2-k_3})}.
\eeqa
De proche en proche, on arrive à
\beqa
\lqn{(\ff^{k_{m-1}}\circ\gg^{k_{m-2}}\circ\dots\circ\gg^{k_2}\circ\ff^{k_1})(j)}
\\[-2ex]
&\hspace{-.5em}=&\hspace{-.5em}\overline{j_{p-k_1-\dots-k_{m-1}-1}\dots j_0}
\\
&&\hspace{-.5em}\overline{(1-j_{p-1})\dots(1-j_{p-k_1})j_{p-k_1-1}\dots j_{p-k_1-k_2}}
\\
&&\hspace{-.5em}\overline{\vphantom{j_p}\dots\vphantom{j_p}}
\\
&&\hspace{-.5em}\overline{(1-j_{p-k_1-\dots-k_{m-2}-1})\dots(1-j_{p-k_1-\dots-k_{m-1}})}
\\[-1.5ex]
\noalign{\newpage}
&\hspace{-.5em}=&\hspace{-.5em}\overline{j_{k_m-1}\dots j_0(1-j_{p-1})\dots(1-j_{p-k_1})}
\\
&&\hspace{-.5em}\overline{j_{p-k_1-1}\dots j_{p-k_1-k_2}}
\\
&&\hspace{-.5em}\overline{\vphantom{j_p}\dots\vphantom{j_p}}
\\
&&\hspace{-.5em}\overline{(1-j_{k_m+k_{m-1}-1})\dots(1-j_{k_m})}
\eeqa
et enfin à
\beqa
\lqn{(\gg^{k_m}\circ\ff^{k_{m-1}}\circ\dots \circ\gg^{k_2}\circ\ff^{k_1})(j)}
\\[-2ex]
&=&\overline{\underset{k_1}{\underbrace{(1-j_{p-1})\dots(1-j_{p-k_1})}}
\underset{k_2}{\underbrace{j_{p-k_1-1}\dots j_{p-k_1-k_2}}}}
\\
&&\overline{\vphantom{j_p}\dots\vphantom{j_p}}
\\
&&\overline{\underset{k_{m-1}}{\underbrace{(1-j_{k_m+k_{m-1}-1})\dots(1-j_{k_m})}}
\underset{k_m}{\underbrace{j_{k_m-1}\dots j_0}}}.
\eeqa
En d'autres termes, les bits de $i$ coïncident avec les bits de $j$ ou
leur complémentaire (le complémentaire d'un bit $\jmath$ étant $1-\jmath$)
selon la règle suivante : de gauche à droite,
\bitem
\item
les $k_1$ premiers bits de $i$ sont les complémentaires de ceux des $k_1$ premiers de $j$,
\item
les $k_2$ bits de $i$ suivants sont identiques aux $k_2$ suivants de $j$,
\item
les $k_3$ bits de $i$ suivants sont les complémentaires des $k_3$ suivants de $j$,
\item
les $k_4$ bits de $i$ suivants sont identiques aux $k_4$ suivants de $j$,
\item
etc.
\eitem
Ce calcul permet d'élaborer un algorithme pour déplacer une carte de numéro donné $j$
vers une position de numéro donné $i$ ($i\neq j$). On décompose $i$ et $j$ en écriture binaire :
$i=\overline{i_{p-1}\dots i_0}$ et $j=\overline{j_{p-1}\dots j_0}$.
Puis on compare les bits de $i$ et $j$ situé à chaque même place
et l'on fait apparaître dans $i$ des blocs de bits successifs identiques à ceux de $j$
et des blocs de bits successifs complémentaires à ceux de $j$.
On décompose ainsi $i$ en <<~blocs de coïncidence~>> et <<~blocs de complémentarité~>> avec ceux de $j$.
Plus précisément, en introduisant la suite des longueurs de ces blocs
\renewcommand{\l}{\lambda}
\beqa
\l_1&=&\min\{k\ge 0:i_k=1-j_k\},
\\
\l_2&=&\min\{k\ge \l_1:i_{k+\l_1}=j_{k+\l_1}\},
\\
\l_3&=&\min\{k\ge \l_2:i_{k+\l_2}=1-j_{k+\l_2}\},
\\
\l_4&=&\min\{k\ge \l_3:i_{k+\l_3}=j_{k+\l_3}\},
\\
&\vdots&
\eeqa
on a, s'il y a $m$ tels blocs ($\l_1+\dots+\l_m=p$ avec
$\l_1,\l_m\ge 0$ et $\l_2,\dots,\l_{m-1}\ge 1$),
\beqa
i&=&\overline{\underset{\l_m}{\underbrace{(1-j_{p-1})\dots(1-j_{p-\l_m})}}}
\\
\noalign{\newpage}
&&\overline{\underset{\l_{m-1}}{\underbrace{j_{p-\l_m-1}\dots j_{p-\l_m-\l_{m-1}}}}}
\\
&&\overline{\vphantom{j_p}\dots\vphantom{j_p}}
\\
&&\overline{\underset{\l_2}{\underbrace{(1-j_{\l_1+\l_2-1})\dots(1-j_{\l_1})}}}
\\
&&\overline{\underset{\l_1}{\underbrace{j_{\l_1-1}\dots j_0}}}.
\eeqa
Les calculs précédents montrent, en choisissant $k_1=\l_m,k_2=\l_{m-1},\dots,k_m=\l_1$, que
$$
(\gg^{\l_1}\circ\ff^{\l_2}\circ\dots \circ\gg^{\l_{m-1}}\circ\ff^{\l_m})(j)=i.
$$
Cela indique que $\l_m$ in-, $\l_{m-1}$ out-, $\dots$, $\l_2$ in- et $\l_1$
out-shuffles mènent la carte \no $j$ à la position \no $i$. En d'autres termes,
le procédé recherché se schématise selon la succession suivante (de gauche à droite) :
$$
\underset{\l_m}{\underbrace{I-I}}\underset{\l_{m-1}}{\underbrace{O-O}}
\dots\underset{\l_2}{\underbrace{I-I}}\underset{\l_1}{\underbrace{O-O}}.
$$
D'un point de vue pratique, on effectue un out-shuffle chaque fois que l'on rencontre,
dans la lecture de gauche à droite des écritures binaires de $i$ et $j$,
une coïncidence de bits et un in-shuffle lorsque l'on rencontre une complémentarité de bits.
Notons que les numéros $i$ et $j$ jouent un rôle symétrique dans cette analyse.
Ainsi, ce processus qui déplace la carte \no $j$ à la place \no $i$ déplace
également la carte \no $i$ à la position \no $j$.

\exam
Considérons le cas d'un jeu de $32$ cartes numérotées $0,1,2,\dots,31$. On souhaite déplacer la carte \no
$19$ vers la position \no $7$. On écrit les nombres $7$ et $19$ en binaire
$7=\overline{00111}$ et $19=\overline{10011}$, on superpose les deux séries
de bits et on compare les paires de bits verticaux. Lorsque l'on a une paire de bits
identiques, on inscrit un out-shuffle ; lorsque l'on a une paire de bits différents,
on inscrit un in-shuffle. Cela donne concrètement :
\newlength{\longueurbit}\settowidth{\longueurbit}{$0$}
\newlength{\hauteurbit}\settoheight{\hauteurbit}{$0$}
$$
\framebox[2\longueurbit]{$\begin{array}{c}0\\1\\-\\I\\[-.4\hauteurbit]\end{array}$}\,
\framebox[2\longueurbit]{$\begin{array}{c}0\\0\\-\\O\\[-.4\hauteurbit]\end{array}$}\,
\framebox[2\longueurbit]{$\begin{array}{c}1\\0\\-\\I\\[-.4\hauteurbit]\end{array}$}\,
\framebox[2\longueurbit]{$\begin{array}{c}1\\1\\-\\O\\[-.4\hauteurbit]\end{array}$}\,
\framebox[2\longueurbit]{$\begin{array}{c}1\\1\\-\\O\\[-.4\hauteurbit]\end{array}$}
$$
L'algorithme pour amener la carte \no $19$ à la place \no $7$ est schématisé
par la succession d'in- et d'out-shuffles $IOIOO$, soit :
$\ff(19)=6$, $\gg(6)=12$, $\ff(12)=25$, $\gg(25)=19$, $\gg(19)=7$.
Globalement,
$$
(\gg^2\circ\ff\circ\gg\circ\ff)(19)=7.
$$
L'algorithme pour amener la carte \no $7$ à la place \no $19$ est le même :
$\ff(7)=15$, $\gg(15)=30$, $\ff(30)\lb =28$, $\gg(28)=25$, $\gg(25)=19$,
qui donne
$$
(\gg^2\circ\ff\circ\gg\circ\ff)(7)=19.
$$
On observe que cette procédure est loin d'être optimale dans le premier cas
puisque seul un out-shuffle suffit à déplacer la carte \no $19$ vers la
position \no $7$ : $\gg(19)=7$...
Néanmoins, elle a l'avantage de fonctionner systématiquement.
D'ailleurs dans le deuxième cas, nous avons vérifié à l'aide de Maple qu'aucune
composition de moins de cinq in-/out-shuffles ne permettait la manipulation
requise ; dans ce cas, l'algorithme se révèle optimal.
Dans~\cite{diaconis-graham}, les auteurs proposent un algorithme minimal
pour exécuter le déplacement souhaité et cet algorithme est valable dans le cas
d'un jeu contenant un nombre quelconque de cartes.

\subsection{\textsf{Déplacement de la carte du dessous du paquet vers une position donnée}}

Examinons le déplacement de la carte du dessous du paquet, \ie\ la carte \no $0$,
vers la position donnée \no $i$. Les calculs précédents donnent immédiatement
\beqa
\lqn{(\gg^{k_m}\circ\ff^{k_{m-1}}\circ\dots\circ\gg^{k_2}\circ\ff^{k_1})(0)}
\\[-2ex]
&=&\overline{\underset{k_1}{\underbrace{1-1}}\underset{k_2}{\underbrace{0-0}}
\dots\underset{k_{m-1}}{\underbrace{1-1}}\underset{k_m}{\underbrace{0-0}}}.
\eeqa
Introduisons la décomposition binaire de $i$ par blocs
$$
i=\overline{\underset{l_m}{\underbrace{1-1}}\underset{l_{m-1}}{\underbrace{0-0}}
\dots\underset{l_2}{\underbrace{1-1}}\underset{l_1}{\underbrace{0-0}}},
$$
où les $l_1,\dots,l_m$ sont les longueurs positives des blocs de bits de $i$
lus de droite à gauche. \'Eventuellement, on posera $l_m=0$ si la décomposition
démarre par un bloc de $0$ et $l_1=0$ si elle finit par un bloc de $1$.
La règle précédente nous enseigne que
$$
(\gg^{l_1}\circ\ff^{l_2}\circ \dots \circ \gg^{l_{m-1}}\circ \ff^{l_m}) (0)=i.
$$
Dans le cas où $l_m=0$, c'est-à-dire dans le cas où la décomposition de $i$ commence
par un bloc de $0$, puisque $\gg(0)=0$ (un out-shuffle n'affecte pas la carte \no $0$),
on peut retirer la manipulation redondante $\gg^{l_{m-1}}$ ci-dessus pour obtenir
$$
(\gg^{l_1}\circ \ff^{l_2}\circ \dots \circ \gg^{l_{m-3}}\circ \ff^{l_{m-2}}) (0)=i.
$$
Ainsi, en effectuant successivement $l_m$ in-, $l_{m-1}$ out-, $\dots$,
$l_2$ in- et $l_1$ out-shuffles, la carte \no $i$ se retrouve au bas du paquet.
C'est le fameux algorithme proposé par Elmsley~\cite{elmsley}. D'un point de vue
pratique, comme cela est mentionné dans~\cite{diaconis-graham},
\cite{diaconis-graham-kantor} et~\cite{elmsley}, on suit le schéma
d'in/out-shuffles dicté par l'écriture binaire de $i$ de gauche à droite
en interprétant un bit~$1$ par un in-shuffle~$I$ et un bit~$0$ par un
out-shuffle~$O$, le premier bloc de $O$ étant omis lorsque $l_m=0$ :
$$
\underset{l_m}{\underbrace{I-I}}\underset{l_{m-1}}{\underbrace{O-O}}
\dots\underset{l_2}{\underbrace{I-I}}\underset{l_1}{\underbrace{O-O}}.
$$

\rem
En fait, la procédure précédemment décrite pour amener la carte \no $0$ à la
position \no $i$ reste valable pour un jeu de $2n$ cartes, $n$ étant un nombre
quelconque. En effet, rappelons que pour $j\in\{0,1,\dots,n-1\}$,
$\ff(j)=2j+1$ et $\gg(j)=2j.$
Soit alors un $j\in\{0,1,\dots,n-1\}$ d'écriture binaire $j=\overline{j_{p-1}\dots j_0}$.
Pour un tel $j$, on a
$$
\ff(j)=\overline{j_{p-1}\dots j_01},\quad
\gg(j)=\overline{j_{p-1}\dots j_00}.
$$
Plus généralement, si
$\overline{j_{p-1}\dots j_0\underset{k}{\underbrace{1-1}}}\le 2n-1$,\lb
(cette condition montrant que
$\overline{j_{p-1}\dots j_0\underset{k-1}{\underbrace{1-1}}}\lb\le n-1$
et que l'utilisation de l'expression de $\ff$ ci-dessus est licite),
alors
$$
\ff^{k}(j)=\overline{j_{p-1}\dots j_0\underset{k}{\underbrace{1-1}}}.
$$
De même, si
$\overline{j_{p-1}\dots j_0\underset{k}{\underbrace{0-0}}}\le 2n-1,$
alors
$$
\gg^{k}(j)=\overline{j_{p-1}\dots j_0\underset{k}{\underbrace{0-0}}}.
$$
Ainsi, pour $l_m$ tel que $\overline{\underset{l_m}{\underbrace{1-1}}}\le 2n-1$, on a
$$
\ff^{l_m}(0)=\overline{\underset{l_m}{\underbrace{1-1}}},
$$
puis, pour $l_{m-1}$ tel que
$\overline{\underset{l_m}{\underbrace{1-1}}\underset{l_{m-1}}{\underbrace{0-0}}}\le 2n-1$, on a
$$
(\gg^{l_{m-1}}\circ\ff^{l_m})(0)=\overline{\underset{l_m}{\underbrace{1-1}}\underset{l_{m-1}}{\underbrace{0-0}}}.
$$
De manière générale, si les nombres $l_1,\dots,l_m$ vérifient
$\overline{\underset{l_m}{\underbrace{1-1}}\underset{l_{m-1}}{\underbrace{0-0}}
\dots\underset{l_2}{\underbrace{1-1}}\underset{l_1}{\underbrace{0-0}}}\le 2n-1$,
alors
$$
(\gg^{l_1}\circ\dots\circ\ff^{l_m})(0)=
\overline{\underset{l_m}{\underbrace{1-1}}\underset{l_{m-1}}{\underbrace{0-0}}
\dots\underset{l_2}{\underbrace{1-1}}\underset{l_1}{\underbrace{0-0}}}.
$$
Cette dernière égalité prouve que la carte \no $0$ peut effectivement atteindre
n'importe quelle position $i\in\{1,2,\dots,2n-1\}.$

\subsection{\textsf{Déplacement de la carte du dessus du paquet vers une position donnée}}

De manière analogue, regardons le déplacement de la carte du dessus du paquet,
\ie\ la carte \no $2^p-1$, vers la position \no $i$.
Puisque $2^p-1$ admet la simple décomposition binaire
$\overline{\underset{p}{\underbrace{1-1}}}$, on a, en inversant l'ordre de $\ff$
et $\gg$ dans les calculs précédents et en se souvenant que $\gg(2^p-1)=2^p-1$,
\beqa
\lqn{(\ff^{k_m}\circ\gg^{k_{m-1}}\circ\dots\circ\ff^{k_2}\circ\gg^{k_1})(2^p-1)}
\\[-3ex]
&=&(\ff^{k_m}\circ\gg^{k_{m-1}}\circ\dots\circ\ff^{k_2})(2^p-1)
\\
&=&\overline{\underset{k_1}{\underbrace{1-1}}\underset{k_2}{\underbrace{0-0}}
\dots\underset{k_{m-1}}{\underbrace{1-1}}\underset{k_m}{\underbrace{0-0}}}.
\eeqa
Donc, pour atteindre la position $i$ avec
$$
i=\overline{\underset{l_m}{\underbrace{1-1}}\underset{l_{m-1}}{\underbrace{0-0}}
\dots\underset{l_2}{\underbrace{1-1}}\underset{l_1}{\underbrace{0-0}}},
$$
on suivra le schéma
$$
\underset{l_{m-1}}{\underbrace{I-I}}\underset{l_{m-2}}{\underbrace{O-O}}
\dots\underset{l_2}{\underbrace{O-O}}\underset{l_1}{\underbrace{I-I}}.
$$

\subsection{\textsf{Déplacement d'une carte donnée vers le dessous du paquet}}

Regardons maintenant le déplacement d'une carte de numéro donné $i$ vers
le dessous du paquet, \ie\ vers la position \no $0$.
Avec la même décomposition binaire de $i$, on a
$$
(\gg^{l_1}\circ\ff^{l_2}\circ \dots \circ \gg^{l_{m-1}}\circ \ff^{l_m}) (i)=0.
$$
Puisque $\gg(0)=0$, on peut omettre la manipulation redondante
$\gg^{l_1}$ ci-dessus pour obtenir
$$
(\ff^{l_2}\circ \gg^{l_3}\circ \dots \circ \gg^{l_{m-1}}\circ \ff^{l_m}) (i)=0.
$$
Ainsi, le schéma suivant déplace la carte \no $i$ au-dessous du paquet :
$$
\underset{l_m}{\underbrace{I-I}}\underset{l_{m-1}}{\underbrace{O-O}}
\dots\underset{l_3}{\underbrace{O-O}}\underset{l_2}{\underbrace{I-I}}.
$$

\subsection{\textsf{Déplacement d'une carte donnée vers le dessus du paquet}}

Enfin, le déplacement de la carte de numéro donné $i$ vers le dessus du paquet,
\ie\ vers la position \no $2^p-1$, est décrit par la relation
$$
(\ff^{l_1}\circ\gg^{l_2}\circ \dots \circ \ff^{l_{m-1}}\circ \gg^{l_m}) (i)=2^p-1.
$$
Ainsi, le schéma suivant, en omettant le dernier bloc de $O$ lorsque $l_1=0$,
déplace la carte \no $i$ au-dessus du paquet :
$$
\underset{l_m}{\underbrace{O-O}}\underset{l_{m-1}}{\underbrace{I-I}}
\dots\underset{l_2}{\underbrace{O-O}}\underset{l_1}{\underbrace{I-I}}.
$$

\section{\textsf{Généralisation : $k$ paquets de $n$ cartes}}\label{section-gene}

On dispose de $k$ paquets de $n$ cartes, donc de $kn$ cartes que l'on numérote
de $1$ à $kn$ ou de $0$ à $kn-1$. Cette situation peut se réaliser avec le
concours de $k$ joueurs installés à une table ronde, ayant chacun un jeu de $n$
cartes en supposant que toutes les cartes sont différentes~: le premier joueur
a un jeu de cartes numérotées de bas en haut $1,2,\dots,n$,
le deuxième a un jeu de cartes numérotées de bas en haut $n+1,n+2,\dots,2n$, etc.,
le $k\ieme$ a un jeu de cartes numérotées de bas en haut $(k-1)n+1,\lb(k-2)n+2,\dots,kn$.
On réalise un in-shuffle ou un out-shuffle en prenant une carte à partir du haut
du paquet à chaque joueur dans un ordre circulaire jusqu'à épuisement des cartes.
On constitue ainsi un nouveau jeu de $kn$ cartes que l'on recoupe en
$k$ paquets de $n$ cartes que l'on redistribue à chaque joueur
et l'on reproduit la manipulation \textit{ad lib.}
Cette généralisation est abordée dans \cite{morris2}.

\subsection{\textsf{In-shuffle}}

Dans le cas d'un in-shuffle généralisé, on commence par prélever la carte de
dessus au premier joueur, puis celle de dessus au deuxième que l'on place sous
la précédente, etc., puis celle de dessus au dernier joueur que l'on place
au-dessous des précédentes. On recommence à partir du premier joueur et
ainsi de suite jusqu'à épuisement des cartes.
On obtient de la sorte un paquet de cartes réparties de bas en haut selon la succession\lb
$(k-1)n+1,(k-2)n+1,\dots,2n+1,n+1,1,$ puis
$(k-1)n+2,(k-2)n+2,\dots,2n+2,n+2,2,$ puis
$(k-1)n+3,(k-2)n+3,\dots,2n+3,n+3,3,$ et ainsi de suite jusqu'à
$kn,(k-1)n,\dots,3n,2n, n$ (voir Fig.~\ref{fig-in-shuffle-generalise}).

\dessininshufflegene

L'in-shuffle généralisé est mathématiquement décrit par la permutation $f$ des
entiers $1,2,\dots,kn$ suivante (voir Fig.~\ref{fig-melanges-generalises}) :
$$
f(i)=\!\left\{\!\!\!\begin{array}{l@{\hspace{.5em}}l}
\dis\frac{i-1}{k}+(k-1)n+1 & \mbox{si \modulo{i}{1}{k},}
\\[2ex]
\dis\frac{i-2}{k}+(k-2)n+1 & \mbox{si \modulo{i}{2}{k},}
\\
\,\,\hphantom{\dis\frac{i-1}{k}}\vdots
\\[2ex]
\dis\frac{i-k+1}{k}+n+1 & \mbox{si \modulo{i}{k\!-\!\!1\!}{\!k},}
\\[1ex]
\dis\frac{i}{k} & \mbox{si \modulo{i}{0}{k}.}
\end{array}\right.
$$
De manière plus condensée, pour $l\in\{0,1,2,\dots,\lb k-1\}$
et $\modtext{i}{l}{k}$,
$$
f(i)=\frac{i-l}{k}+(k-l)n+1.
$$

\dessininshufflepermutations

La réciproque de $f$ est donnée par
$$
f^{-1}(j)=\!\left\{\!\!\!\begin{array}{l@{\hspace{-2em}}l}
kj & \mbox{si $1\le j\le n$,}
\\[1ex]
kj-(kn+1) & \mbox{si $n+1\le j\le 2n$,}
\\[1ex]
kj-2(kn+1) & \mbox{si $2n+1\le j\le 3n$,}
\\
\,\,\hphantom{kj-}\vdots
\\[1ex]
kj-(k-1)(kn+1)
\\
& \hspace{-1.5em}\mbox{si $(k-1)n+1\le j\le kn$.}
\end{array}\right.
$$
La forme générique est, pour $l\in\{0,1,2,\dots,k-1\}$
et $ln+1\le j\le (l+1)n$,
$$
f^{-1}(j)= kj-l(kn+1).
$$
On note en particulier la congruence
$$
\modulo{f^{-1}(j)}{kj}{(kn+1)}.
$$
Cette observation conduit à la formulation suivante de la période de $f$.
\bth{
La période de $f$ est l'ordre de $k$ modulo $(kn+1)$, c'est-à-dire le premier
entier \mbox{$r\ge 1$} tel que $\modtext{k^r}{1}{(kn+1)}$.
}

\bco{Si $n$ est de la forme $k^{p-1}$ pour un $p\ge 1$, alors la période
de $f$ est $2p$.
}
\dem
Lorsque $n=k^{p-1}$, alors $kn+1\linebreak =k^p+1$ et $\modtext{k^p}{-1}{(kn+1)}$.
On a donc $\modtext{k^{2p}}{1}{(kn+1)}$ qui prouve que $f^{2p}=id$.
Par ailleurs, l'orbite de la carte \no $1$ sous l'action de la permutation $f$ est
\beqa
\cO(1)&=&\{1,k,k^2,\dots,k^{p-1},k^p-k^{p-1},\\
&&
k^p-k^{p-2},\dots,k^p-k,k^p-1\}
\eeqa
qui est de cardinal $2p$. C'est la période de $f$.
\fin

\subsection{\textsf{Out-shuffle}}

Dans le cas d'un out-shuffle généralisé, on commence par prélever la carte
de dessus au dernier joueur, puis celle de dessus à l'avant-dernier que l'on place
sous la précédente, etc., puis celle de dessus au premier joueur. On
recommence à partir du dernier joueur et ainsi de suite jusqu'à épuisement des cartes.
On obtient alors la répartition des cartes de bas en haut selon la succession
$0,n,2n,\dots, (k-2)n, (k-1)n,$ puis
$1,n+1,2n+1,\dots, (k-2)n+1,(k-1)n+1,$ puis
$2,n+2, 2n+2,\dots,(k-2)n+2,(k-1)n+2,$\lb et ainsi de suite jusqu'à
$n-1,2n-1,3n-1,\dots,\lb (k-1)n-1, kn-1$ (voir Fig.~\ref{fig-out-shuffle-generalise}).

\dessinoutshufflegene

L'out-shuffle généralisé est modélisé par la permutation $\tg$ des entiers
$0,1,2,\dots,kn-1$ suivante (voir Fig.~\ref{fig-melanges-generalises}) :
$$
\tg(i)=\!\left\{\!\!\!\begin{array}{l@{\hspace{.5em}}l}
\dis\frac{i}{k} & \mbox{si \modulo{i}{0}{k},}
\\[2ex]
\dis\frac{i-1}{k}+n & \mbox{si \modulo{i}{1}{k},}
\\[2ex]
\dis\frac{i-2}{k}+2n & \mbox{si \modulo{i}{2}{k},}
\\
\,\,\hphantom{\dis\frac{i-1}{k}}\vdots
\\[2ex]
\dis\frac{i-k+1}{k}+(k-1)n & \mbox{si \modulo{i}{k\!-\!1\!}{\!k}.}
\end{array}\right.
$$
La réciproque de $\tg$ est donnée par
$$
\tg^{-1}(j)=\!\left\{\!\!\!\begin{array}{l@{\hspace{-2em}}l}
kj & \mbox{si $0\le j\le n-1$,}
\\[1ex]
kj-(kn-1) & \mbox{si $n\le j\le 2n-1$,}
\\[1ex]
kj-2(kn-1) & \mbox{si $2n\le j\le 3n-1$,}
\\
\,\,\hphantom{kj-}\vdots
\\[1ex]
kj-(k-1)(kn-1)
\\
&\hspace{-1.3em} \mbox{si $(k-1)n\le j\le kn-1$.}
\end{array}\right.
$$
On a en particulier la congruence
$$
\modulo{\tg^{-1}(j)}{kj}{(kn-1)}
$$
qui conduit à la formulation suivante de la période de $\tg$.

\bth{
La période de $\tg$ est l'ordre de $k$ modulo $(kn-1)$, c'est-à-dire le premier
entier \mbox{$s\ge 1$} tel que \modtext{k^s}{1}{(kn-1)}.
}

\bco{Si $n$ est de la forme $k^{p-1}$ pour un $p\ge 1$, alors la période
de $\tg$ est $p$.
}
\dem
Lorsque $n=k^{p-1}$, on a $kn-1\linebreak =k^p-1$ et alors $\modtext{k^p}{1}{(kn+1)}$.
On a donc $\tg^p=id$.
L'orbite de la carte \no $1$ sous l'action de la permutation $\tg$ est
simplement
$$
\tcO(1)=\{1,k,k^2,\dots,k^{p-1}\}
$$
qui est de cardinal $p$. C'est la période de $\tg$.
\fin

\noindent\textsc{Remerciements.}
J'adresse mes sincères remerciements à deux de mes élèves, Matthieu Bacconnier
et Anthony Tschirhard (INSA de Lyon, 51\ieme\ promotion), le premier pour son
aide relative aux calculs numériques présentés dans la section~\ref{section-num},
le second pour m'avoir soumis ces problèmes de mélange qui auront abouti au
présent travail. D'autres remerciements s'adressent à Philippe Biane pour m'avoir
communiqué certaines références sur le sujet.

\def\refname{Références}

\end{document}